\documentclass[reqno,a4paper]{amsart}

\usepackage{style}

\makeindex

\title{Localisations of rings of weak global dimension at most one}

\author{Jan {\v{S}}{\v{t}}ov\'{\i}\v{c}ek}
\address{Department of Algebra, Faculty of Mathematics and Physics, Charles University, Sokolovsk\'{a} 83, 186~75 Prague, Czech Republic}
\email{stovicek@karlin.mff.cuni.cz}

\author{Florian Tecklenburg}
\address{Fakultat f\"ur Mathematik, Universit\"at Bielefeld, D-33501 Bielefeld, Germany}
\email{ftecklen@math.uni-bielefeld.de}

\date{\today}
\begin{document}

\maketitle
\begin{abstract}
%
We characterise the commutative rings of weak global dimension at most one whose derived categories can be stratified by the Balmer--Favi support theory.
For such rings, this yields a bijective correspondence between the class of localising subcategories of the derived category, and the subsets of the Zariski spectrum.
For a large class of rings of weak global dimension at most one including all valuation domains, this occurs if and only if the ring is locally strongly discrete, if and only if the Telescope conjecture holds, and if and only if $\MOD_R$ has Gabriel dimension.
%

However, we also give an example of a valuation domain where all these conditions fail, and yet the residue fields generate its derived category and we can classify all the localising subcategories of $\sfD(R)$. Additionally, in the same example, we show that the Balchin--Stevenson smashing support cannot stratify $\sfD(R)$ either.
\end{abstract}
\tableofcontents

\section{Introduction}
Given a commutative ring $R$, we are interested in understanding the lattice of localising subcategories of the unbounded derived category over $R$: $\Loc(\sfD(R))$. There are only a few classes of examples known in which this lattice has been fully classified. For instance, take $R$ to be a noetherian ring, or a semi-artinian absolutely flat ring. For these, the classification manifests itself as an isomorphism of lattices: $\Loc(\sfD(R))\cong\cP(\Spec(R))$, where $\cP(-)$ denotes the power-set functor.
Moreover, such an isomorphism is obtained using a standard support theory. Derived categories admitting this kind of classification are called stratified; this terminology goes back to~\cite{bik11}.
In several other cases, only partial results are known \cite{neemanchromtower,stevenson_absolutelyflatrings,isolated}.

Our work extends the list of examples where stratification can be obtained. 
In the sense of the abstract stratification theory within tensor-triangular geometry using Rickard idempotents \cite{ric97,balmerfavi_idempotents,bhs23}, we give a complete characterisation of when the derived category is stratified for classes of rings of weak global dimension at most one. We refer to a key class of these rings as follows.
\begin{definition*}
    A commutative ring $R$ of weak global dimension at most one is called {\it strongly discrete} if for any prime $\fp$ the localisation $R_\fp$ is a strongly discrete valuation domain.
    I.e., $R_\fp$ is a domain, its ideals are totally ordered by inclusion, and only the trivial ideals $(0)$ and $R_\fp$ are idempotent.

In the case that $R$ is not local (so not a valuation domain itself), we require more than $R$ being strongly discrete in our results. To this end, we call a flat ring epimorphism $f\colon R\to S\neq0$ {\it sharp} (see \ \Cref{lem:sharpepi}) if there exists
a finitely generated ideal $J\subseteq S$ such that $S/J$
is a local ring. With this, we can formulate our first main result.
\end{definition*}
\begin{theorem*}[{\Cref{thm:wgld1strat}}]
    Let $R$ be a ring of weak global dimension at most one. Then the following are equivalent:
    \begin{enumerate}
        \item $R$ is strongly discrete and each non-zero flat epimorphism $R\to S$ is sharp.
        \item $\sfD(R)$ is stratified by $\Spc(\sfD(R)^\omega)$ in the sense of~\cite{bhs23}.
    \end{enumerate}
\end{theorem*}
Furthermore, we provide an example of a ring $R$ of weak global dimension one (even of a valuation domain) that is not strongly discrete, and yet we have an abstract isomorphism of lattices $\Loc(\sfD(R))\cong\cP(\Spec(R))$. 
For this ring $R$, $\Spec(R)$ is not a weakly noetherian space in the sense of~\cite{bhs23}, the Gabriel dimension of $\MOD_R$ is undefined, the Telescope Conjecture for $\sfD(R)$ fails, yet the residue fields generate $\sfD(R)$ as a localising category. 
We also show that we cannot save the day by the finer smashing support theory in the sense of~\cite{bigspectra}, as the smashing support does not stratify $\sfD(R)$ either. This is all \Cref{ex:ce}. 

Finally, our second main result resolves the classification problem of $\Loc(\sfD(R))$ via stratification by the Zariski spectrum for a class of coherent rings $R$ of weak global dimension at most one%
\footnote{When writing up our results, we learned that Barthel--Heard-Sanders-Zou provided an independent proof of the stratification of strongly discrete valuation domains in~\cite{BHSZ}.}, e.g., Prüfer domains.
We also put stratification problem into context with Gabriel dimension of the modules categories, partially extending results of Balchin--Gómez--Stevenson~\cite{isolated}.

\begin{theorem*}[{\Cref{thm:kg=stratforvd}, \Cref{cor:kgandwgld}}]
    Let $R$ be a coherent ring of weak global dimension at most one such that one of the following conditions holds:
    \begin{tasks}[style=enumerate, label-align=right ,column-sep = {-7mm},item-indent = {7mm},label-offset = {0.8em} ](1)
        \task*[(a)] any non-zero flat epimorphism is sharp, or
        \task*[(b)] $R$ is a valuation domain.
    \end{tasks}
    Then the following are equivalent:
\begin{tasks}[style=enumerate, label-align=right ,column-sep = {-7mm},item-indent = {7mm},label-offset = {0.8em} ](2)
        \task*(1)[(i)] $R$ is strongly discrete.
        \task*[(ii)] $\sfD(R)$ satisfies the Telescope Conjecture.
        \task*(1)[(iii)] $\sfD(R)$ is stratified by $\Spc(\sfD(R)^\omega)$.
        \task*[(iv)] $\MOD_R$ has Gabriel dimension.
    \end{tasks}
\end{theorem*}

All considered rings are commutative and unital unless stated otherwise. If $R$ is a ring and $\fp\in\Spec(R)$, we denote the residue field $R_\fp/\fp R_\fp$ at $\fp$ by $k(\fp)$. If $R$ is a domain, $Q$ stands for the fraction field (= the residue field at the zero ideal).

 \subsection*{Acknowledgements \& AI disclosure} 
The first-named author was supported by the Czech Science Foundation grant GA\v{C}R~26-22734S.
The second-named author was supported by the Deutsche Forschungsgemeinschaft (Project-ID 491392403 -- TRR 358).

AI usage: Used solely for finding suitable references to the existing literature, spelling improvements and correction of typos.

\section{Preliminaries}

\subsection{Commutative domains and ring epimorphisms}

As alluded to, our main results are about rings of weak global dimension at most one. Here we will recall the necessary concepts about non-noetherian commutative rings and in particular domains, following the terminology of~\cite{fs01}.

\begin{definition}[{\cite[\S II.1]{fs01}}] \label{def:VD}
By a \emph{valuation domain}, we mean a commutative domain whose ideals are ordered by inclusion.
\end{definition}

Standard but from the point of view here uninteresting, examples are discrete valuation rings (or discrete rank $1$ valuation domains in the terminology of \cite{fs01}) such as the ring of $p$-adic integers $\Z_p$ for a prime $p$ or the power series ring $k[[x]]$ for a field $k$.

\begin{definition}[{\cite[\S II.3]{fs01}}] \label{def:value-group}
Given a valuation domain $R$ with the fraction field $Q$, the quotient group $\Gamma(R) = Q^\times/R^\times$ is called the \emph{value group} of $R$. Here, $R^\times$ and $Q^\times$ are the multiplicative groups of units of $R$ and $Q$, respectively.
\end{definition}

In fact, the value group is a totally ordered abelian group, where we define the ordering by putting $aR^\times\le bR^\times$ for $a,b\in Q^\times$ if $ba^{-1}\in R$, or equivalently when $bR\subseteq aR$ as $R$-submodules of $Q$.
Moreover, all totally ordered abelian groups occur as value groups of valuation domains by a result of Krull~\cite{kru32} (see also~\cite[Theorem~II.3.8]{fs01}).

We will prove best results for strongly discrete valuation domains. We recall the corresponding definitions from~\cite[\S II.8]{fs01}.

\begin{definition} \label{def:discrete-VD} Let $R$ be a valuation domain. We call $R$
\begin{enumerate}
\item \emph{discrete} if for any two adjacent prime ideals $\fp\subset\fq$, the factor ring $R_\fq/\fp$ is a discrete rank $1$ valuation domain (notice that $\fp R_\fq = \fp$ as subsets of $R_\fq$ since $\fp$ is canonically an $R_\fq$-module by~\cite[Lemma~4.3(iv)]{fs01});
\item \emph{strongly discrete} if no non-zero prime of $R$ is idempotent.
\end{enumerate}
\end{definition}

The relation between the two just defined notions is made precise by the following.

\begin{theorem}[{\cite[\S II Theorem 8.3]{fs01}}]\label{thm:propertiesStrong}
For a valuation domain $R$, the following are equivalent:
\begin{enumerate}
    \item $R$ is strongly discrete (i.e.\ no non-zero prime ideal in $R$ is idempotent).
    \item For each prime ideal $\fp$, the unique maximal ideal $\fp R_\fp$ is principal.
    \item $R$ is discrete and the prime ideals of $R$ satisfy the ascending chain condition.
\end{enumerate}
\end{theorem}

Valuation domains are, of course, local. For non-local commutative rings, we have the following result very much akin to the characterization of Dedekind domains in the noetherian case.

\begin{lemma}[{\cite[Lemma 4.2.6]{Glaz:1989}}]\label{lem:wgld1}
    A commutative ring $R$ has weak global dimension at most one if and only if $R_\fp$ is a valuation domain for each prime ideal~$\fp$, if and only if $R_\fm$ is a valuation domain for each maximal ideal~$\fm$.
\end{lemma}

In the case of domains, the following terminology is well established.

\begin{definition}[{\cite[\S III.1]{fs01}}] \label{def:Pruefer}
A \emph{Pr\"ufer domain} is a commutative domain such that $R_\fm$ is a valuation domain for any maximal ideal $\fm$.
\end{definition}

Except for domains themselves, we will also need to understand certain class of ring homomorphisms between them.

\begin{definition}
A ring homomorphism $f\colon R\to S$ between (a priori) non-commutative rings is
\begin{enumerate}
\item a \emph{ring epimorphism} if it is an epimorphism in the category of rings (i.e.\ whenever $g,h\colon S\to T$ are homomorphism of rings such that $gf=hf$, then $g=h$);
\item a \emph{homological epimorphism} if it is a ring epimorphism such that $\Tor^R_n(S,S)=0$ for all $n>0$;
\item a \emph{flat epimorphism} if it is a ring epimorphism such that $S$ is flat as a left $R$-module.
\end{enumerate}

Moreover, two ring epimorphisms $f_i\colon R\to S_i$, $i=1,2$ are called \emph{equivalent} if there exists a (necessarily unique) ring isomorphism $g\colon S_1\to S_2$ such that $gf_1=f_2$
\end{definition}

A ring homomorphism $f\colon R\to S$ is a ring epimorphism if and only if the multiplication in $S$ induces an isomorphism $S\otimes_RS\to S$; see \cite[Proposition~1.1]{Silver:1967}. Each flat epimorphism is clearly a homological epimorphism but not vice versa---an explicit example was given in~\cite{keller94}.
Moreover, if $R$ is commutative and $f\colon R\to S$ is a ring epimorphism, then $S$ is commutative by \cite[Corollary~1.2]{Silver:1967} and the induced map of Zariski spectra $\Spec(f)\colon\Spec(S)\to\Spec(R)$ is injective by~\cite[Proposition~IV.1.5]{Lazard:1969}.

A typical example of a flat epimorphism of commutative rings is a localisation $R\to \Sigma^{-1}R$ at a multiplicative subset $\Sigma\subset R$, and flat epimorphisms originating from a valuation domain are all of this form by~\cite[Proposition 5.4]{gldm1}.
If $f$ is a flat epimorphism, then $\Spec(f)$ is even a homeomorphism onto its image and $f$ is determined by the image of $\Spec(f)$ up to equivalence~\cite[Corollaire~IV.2.2]{Lazard:1969}.
The latter is, however, not necessarily the case for homological epimorphisms.
Indeed, the example from~\cite{keller94} is a valuation domain $R$ whose spectrum has two points: the zero ideal and the maximal ideal $\fm$. Then the canonical map $f\colon R\to Q\times k(\fm)$ is a homological epimorphism, the map $\Spec(f)$ is a continuous bijection, but it is not a homeomorphism since $\Spec(Q\times k(\fm))$ is discrete while $\Spec(R)$ is not.
See also~\cite[Example~5.24]{gldm1}.

\subsection{Hochster duality}
\label{subsec:Hochster}

Given a topological space $X$, we denote by $\Omega(X)$ the collection of all open subsets of $X$ and by $\Omega(X)^\omega\subseteq\Omega(X)$ the set of all quasi-compact open sets of $X$.

It is a classical result by Hochster~\cite{hochster} that a topological space $X$ is homeomorphic to one of the form $\Spec(R)$ for a commutative ring $R$ if and only if it is a spectral space in the following sense:

\begin{definition} \label{def:spectral}
A space $X$ is called \emph{spectral} if it satisfies all of the following:
\begin{enumerate}
\item $X$ is $T_0$ and quasi-compact;
\item the quasi-compact open subsets of $X$ are closed under finite intersection and form an open basis;
\item every non-empty irreducible closed subset has a generic point.
\end{enumerate}
\end{definition}

If $X$ is a spectral space, we can define a \emph{Hochster dual space} $X^\vee$ with the same set of points, but whose quasi-compact open sets are precisely complements of the quasi-compact open subsets of $X$. Then $X^\vee$ is again a spectral space and $(X^{\vee})^\vee=X$.
This is still due to Hochster~\cite{hochster} and details can be found in~\cite[Ch.~1]{dst19}.
The open sets of $X^\vee$ are often called the \emph{Thomason subsets} of $X$, paying tribute to Thomason's seminal paper~\cite{Thomason:1997} classifying the thick subcategories of the perfect derived category of $R$ in terms of Thomason subsets of the Zariski spectrum.

By unraveling the definition, Thomason subsets of $\Spec(R)$ are precisely those of the form $\bigcup_{\lambda\in\Lambda}V(I_\lambda)$, where $\{ I_\lambda\subseteq R\mid \lambda\in\Lambda\}$ is a collection of finitely generated ideals of $R$. Thomason subsets are always closed under specialization in $\Spec(R)$, but not vice versa in general.

In connection with classification of localising subcategories of $\sfD(R)$, we will need the following definition. 

\begin{definition} \cite[Definitions~2.1 and~2.3]{bhs23}
A point $x$ in a spectral space $X$ is said to be \emph{weakly visible} if the singleton $\{x\}$ is the intersection of a Thomason subset and the complement of a Thomason subset.
A spectral space $X$ is \emph{weakly noetherian} if each $x\in X$ is weakly visible.
\end{definition}

In other words, $X$ is weakly noetherian if and only if each point $x\in X^\vee$ is locally closed. The latter is a separation axiom between $T_0$ and $T_1$; this connection has been observed in~\cite{bigspectra}.

\begin{definition}\label{def:T_D}
A topological space $Y$ satisfies the \emph{$T_D$ separation axiom}, hence called a \emph{$T_D$-space}, if each singleton $\{x\}\subseteq Y$ is locally closed in $Y$.
\end{definition}

\subsection{Tensor-triangular geometry}

Given any essentially small triangulated category $\cK$ with a symmetric monoidal structure $(\cK,\otimes,\unit)$ such that $\otimes$ is exact in each variable, Balmer defined in \cite{balmer_spectrum} its spectrum $\Spc(\cK)$, nowadays called the \emph{Balmer spectrum}.
As he showed in \cite{balmer-ICM}, this allows one to uniformly treat a large class of examples from algebraic geometry, homotopy theory, representation theory, or motivic theory using geometric intuition.

More in detail, the points of $\Spc(\cK)$ are prime thick tensor ideals of $\cK$. Here $\cI\subseteq\cK$ is a \emph{thick tensor ideal} if $\cI$ is a thick subcategory of $\cK$ such that $x\otimes y\in\cI$ whenever $x\in\cK$ and $y\in\cI$. A thick tensor ideal $\cP\subset\cK$ is called prime if $x\otimes y\in\cP$ implies $x\in\cP$ or $y\in\cP$ for each pair $x,y\in\cK$. The Balmer spectrum carries the \emph{Zariski topology} that makes it a spectral space and whose quasi-compact open sets are precisely those of the form
\[ D(x) = \{ \sfP\in\Spc(\cK) \mid x\in P \}, \]
where $x\in\cK$. The corresponding Thomason sets are of the form
\[ \{ \sfP\in\Spc(\cK) \mid x\not\in\cP \text{ for some } x\in\cS \}, \]
where $\cS$ is a set of objects of $\cK$.

Our prime example of $\cK$ is as follows: If $R$ is a ring, the unbounded derived category $\sfD(R)$ is compactly generated and carries a monoidal structure $(\sfD(R),\dert_R,R)$ such that $\dert_R$ is exact in both variables. This is a compactly generated triangulated category and the full subcategory of compact objects $\cK=\sfD(R)^\omega$ (which are known to be perfect complexes) is an essentially small category as above with the restricted structure. The key result in that context is as follows.

\begin{theorem}\label{thm:Spec-to-Spc}
Let $R$ be a ring. Then there is a homeomorphism $\Spec(R)\simeq\Spc(\sfD(R)^\omega)$ that sends each $\fp\in\Spec(R)$ to the prime thick tensor ideal
\[ \sfP = \{ x\in\sfD(R)^\omega \mid x_\fp \simeq 0 \}. \]
\end{theorem}

\begin{proof}
This goes back to \cite[Corollary~5.6]{balmer_spectrum} for $R$ noetherian, and in general it was shown in~\cite[Theorem~9.5]{Buan-Krause-Solberg:support-var}; see also \cite{balmer-ICM}.
\end{proof}

Next, we recall fragments of the results from \cite{balmerfavi_idempotents,  bhs23} which allow us to classify all localising subcategories or localising tensor ideals in some cases.
To this end, let $(\cT,\otimes,\unit)$ be a compactly generated triangulated category such that $\otimes$ is exact and preserves coproducts in each variable, and denote the full subcategory of compact objects by $\cT^\omega$.
We call a subcategory $\sfL\subseteq\cT$ \emph{localising} if it is closed under arbitrary direct sums, and we call $\sfL$ a \emph{localising tensor ideal} if $\sfL$ is a localising subcategory such that $x\otimes y\in\sfL$ whenever $x\in\cT$ and $y\in\sfL$.
We will denote the smallest localising (resp., tensor) subcategory containing a class of objects $\cS\subseteq\cT$ by $\loc(\cS)$ (resp., $\loc^\otimes(\cS)$). Note that if $\loc(\unit)=\cT$, then every localising subcategory of $\cT$ is a localising tensor ideal by~\cite[Lemma~1.4.6]{HPS:1997}. In particular, given a ring $R$, this holds for $\cT=\sfD(R)$.

By the Brown representability theorem \cite[Theorem~3.1]{neemanBrownRep}, the monoidal structure on $\cT$ is closed under the above assumptions. That is, for each $y,z\in\cT$ there is a hom-object $\RHom(y,z)\in\cT$ such that for each $x,y,z\in\cT$ there is a natural isomorphism
\[ \Hom_{\cT}(x\otimes y,z)\simeq\Hom_{\cT}(x,\RHom(y,z)). \]
This allows us to define \emph{dualisable objects} of $\cT$ as those $x\in\cT$ such that the canonical homomorphism $\RHom(x,\unit)\otimes y\to \RHom(x,y)$ is an isomorphism for each $y\in\cT$. We then call $\cT$ \emph{rigidly-compactly generated} if the classes of compact and dualisable objects in $\cT$ coincide. In that case, $\cT^\omega$ is automatically an essentially small triangulated category, the monoidal structure restricts to $\cT^\omega$, and it is well-known that $\cT=\sfD(R)$ is rigidly-compactly generated. We refer to \cite{balmerfavi_idempotents,  bhs23} for more details.

Given a rigidly-compactly generated category $\cT$ and $\sfP\in\Spc(\cT^\omega)$ weakly visible, Balmer and Favi define in~\cite[Definition~7.9]{balmerfavi_idempotents} a certain tensor idempotent object $g(\sfP)\in\cT$ (i.e.\ one has that $g(\sfP)\otimes g(\sfP)\simeq g(\sfP)$).
Paying tribute to a similar idea of Rickard~\cite{ric97} this time, we will call such $g(\sfP)$ the \emph{Rickard idempotent} associated with $\sfP$. We will not reproduce the general definition here, but show instead how to construct such an object directly for $\cT=\sfD(R)$.
Note also that  in that case $\sfP$ corresponds to a unique prime $\fp\in\Spec(R)$ by \Cref{thm:Spec-to-Spc}, and $\fp$ being weakly visible in $\Spec(R)$ precisely means that there is a finitely generated ideal $I\subseteq R$ such that $V(IR_\fp)=\{\fp R_\fp\}$ in $\Spec(R_\fp)$.

\begin{proposition}\label{prop:g(P)-in-D(R)}
Suppose that $I=(x_1,\dots,x_n)\subseteq R$ is a finitely generated ideal such that $V(IR_\fp)=\{\fp R_\fp\}$ in $\Spec(R_\fp)$. If $\sfP\in\Spc(\sfD(R)^\omega)$ is the point corresponding to $\fp$ in the sense of \Cref{thm:Spec-to-Spc},
then
\[
g(\sfP) \simeq \bigotimes_{i=1}^n\left(R_\fp\to R_\fp\left[\frac{1}{x_i}\right]\right) \in \sfD(R),
\]
where each $R_\fp\to R_\fp\big[\frac{1}{x_i}\big]$ is viewed as a two-term complex with $R_\fp$ sitting in cohomological degree zero.
\end{proposition}

\begin{proof}
Upon identifying $\Spc(\sfD(R)^\omega)$ with $\Spec(R)$, the two Thomason sets of $\Spec(R)$ from \cite[Lemma~7.8]{balmerfavi_idempotents}, witnessing that $\fp$ is weakly visible, are $U=V(I)$ and $W=\Spec(R)\setminus\Spec(R_\fp)$ (in the sense that $\{\fp\}$ equals the intersection of $U$ and the complement of $W$).
Then $g(\sfP)$ is introduced in \cite[Definition~7.9]{balmerfavi_idempotents} as $e(U)\dert_R f(W)$, where $e(U)$ is the so-called left idempotent corresponding to $U$ and $f(W)$ is the so-called right idempotent corresponding to $W$.
It follows from~\cite[Proposition~6.10]{DG:2002} that 
\[
e(\sfP) \simeq \bigotimes_{i=1}^n\left(R\to R\left[\frac{1}{x_i}\right]\right) \in \sfD(R).
\]
Finally, one observes that the localisation functor corresponding to $W$ is none other than the usual localisation at $\fp$, and so $f(W)\simeq R_\fp$ in $\sfD(R)$. It follows that $g(\sfP)=e(U)\dert_R f(W)$ is the complex from the statement.
\end{proof}

If now $\cT$ is rigidly-compactly generated and $\Spc(\cT^\omega)$ is weakly noetherian, we can define a support of each object $x\in\cT$.

\begin{definition} \label{def:support}
The \emph{(Balmer-Favi) support} of $x\in\cT$ is defined as
\[
\supp(x) = \{\sfP\in\Spc(\cT^\omega) \mid g(\sfP)\otimes x\not\simeq 0 \}.
\]
\end{definition}

Finally, we recall a key definition and result for us from~\cite{bhs23}.

\begin{definition}\label{def:stratified}
Given a rigidly-compactly generated triangulated category $\cT$ such that $\Spc(\cT^\omega)$ is weakly noetherian, we call $\cT$ \emph{stratified} if it satisfied the following two conditions:

\begin{itemize}
\item \emph{The local-to-global principle}:
For any object $x\in \cT$ we have 
\[
x \in \loc^\otimes( g(\sfP)\otimes x \mid \sfP\in\Spc(\cT^\omega) ).
\]
\item \emph{The minimality condition}:
For each $\sfP \in \Spc(\cT^\omega)$, the localising tensor ideal
$g(\sfP)\otimes\cT$ is a minimal localising tensor ideal in the sense that it
is generated by any non-zero object it contains.
\end{itemize}
\end{definition}

\begin{theorem}[{\cite[Theorem 4.1]{bhs23}}] \label{thm:stratification-abstract}
Let $\cT^\omega$ be a rigidly-compactly generated triangulated category such that $\Spc(\cT^\omega)$ is weakly noetherian. Then $\cT$ is stratified if and only if the map
\[
\{ \text{localising tensor ideals of } \cT \} \to
\{ \text{subsets of } \Spc(\cT^\omega) \}
\]
given by $\sfL \mapsto \bigcup_{x\in\sfL}\supp(x)$ is a bijection, with inverse $U\mapsto\{x\in\cT\mid\supp(X)\subseteq U\}$.
\end{theorem}

In other words, for stratified categories, one can classify all localising subcategories in terms of the Balmer-Favi support.
If, moreover, $\cT=\sfD(R)$, then localising tensor ideals coincide with localising subcategories by the discussion above.
It is a classical result that $\sfD(R)$ is stratified if $R$ is noetherian~\cite[Theorem~2.8]{neemanchromtower}, but it is also well-known that stratification for $\sfD(R)$ may badly fail if $R$ is not noetherian, \cite{neemanOddball}.
In the rest of the paper, we will characterise when $\sfD(R)$ is stratified for valuation domains and, more generally, for a large class of rings of weak global dimension bounded by one.

\section{Stratification for valuation domains}

Our aim in the section is to prove that for a valuation domain $R$, the category $\sfD(R)$ is stratified if and only if $R$ is strongly discrete.
We begin by proving the local-to-global principle for strongly discrete valuation domains.

\begin{lemma} \label{lem:T_D-for-VD}
If $R$ is a valuation domain, then the following conditions are equivalent:
\begin{enumerate}
\item The prime ideals of $R$ satisfy the ascending chain condition.
\item The Zariski spectrum $\Spec(R)$ is a weakly noetherian topological space (equivalently: $\Spec(R)^\vee$ is a $T_D$-space).
\item The Zariski spectrum is a noetherian topological space (i.e.\ the open subsets satisfy the ascending chain condition).
\end{enumerate}
\end{lemma}
\begin{proof}
Every noetherian topological space is weakly noetherian by~\cite[Remark~2.2]{bhs23}.
If $\Spec(R)$ is weakly noetherian, then the primes of $R$ satisfy the ascending chain condition by \cite[Proposition~4.13]{cosupport}.

Finally, suppose that the primes of $R$ satisfy the ascending chain condition.
Proving that $\Spec(R)$ is noetherian is the same as proving that each non-empty Zariski open subset $\varnothing\ne U\subseteq\Spec(R)$ is quasi-compact.
Now, if $U\subseteq\Spec(R)$ is Zariski open, it is a down-set with respect to the inclusion of primes, and by the ascending chain condition there is a maximal prime ideal $\fp_U\in U$ with respect to inclusion.
In particular, we must have $U=\{\fq\in\Spec(R) \mid \fq\subseteq\fp_U\}$.
Thus, if $U=\bigcup_{\lambda\in\Lambda} W_\lambda$ is an open cover of $U$, we must have $\fp_U\in W_\lambda$ for some $\lambda\in\Lambda$ and, consequently, $U=W_\lambda$.
\end{proof}

\begin{lemma}\label{prop:l2gForStrongDiscrete}
If $R$ is strongly discrete, the local-to-global principle holds for $\sfD(R)$.
\end{lemma}
\begin{proof}
If $R$ is strongly discrete, then $\Spec(R)$ is a noetherian space by \Cref{thm:propertiesStrong} and \Cref{lem:T_D-for-VD}.
Hence also $\Spc(\sfD(R)^\omega)$ is noetherian by~\Cref{thm:Spec-to-Spc}, and $\sfD(R)$ satisfies the local-to-global principle by \cite[Theorem 6.9]{tensoractions}.
\end{proof}

After proving the local-to-global principle, our next goal is to prove minimality of the localising subcategories generated by the Rickard idempotents $g(\sfP)$ (recall \Cref{prop:g(P)-in-D(R)}). To this end, let us fix some notation. The set of primes $(\Spec(R),\subseteq)$ is totally ordered by inclusion, is non-empty and satisfies the ascending chain condition by~\Cref{thm:propertiesStrong}.
If we denote by $(I,\le)$ a totally ordered set opposite to $(\Spec(R),\subseteq)$, we can use the anti-isomorphism $I\simeq\Spec(R)^{\op}$ to index primes; then $\fp_j\subset\fp_i$ whenever $i<j$ in $I$.
If we denote $0=\min I$ and $\infty=\max I$, then $\fp_0$ is the maximal ideal and $\fp_\infty=(0)$.
Similarly, using the inclusion-reversing bijection $\Spec(R)\cong \Spc(\sfD(R)^\omega)$ from \Cref{thm:Spec-to-Spc}, we index elements of $\Spc(\sfD(R)^\omega)$ by $I$ in such a way that $\sfP_0=\{0\}$, $\sfP_\infty$ consists of all perfect complexes with torsion cohomology, and $\sfP_i\subset\sfP_j$ whenever $i<j$ in $I$.

\begin{lemma}\label{lem:BF}
For a strongly discrete valuation domain $R$, we have $g(\sfP_\infty)=Q$, and for each $i<\infty$, $
g(\sfP_i)\simeq(R_{\fp_i}\to R_{\fp_{i+1}})$.
\end{lemma}
\begin{proof}
If $i=\infty$, then $\fp_\infty=0$ and we can put $n=1$ and $x_1=0$ in \Cref{prop:g(P)-in-D(R)}, giving $g(\sfP_\infty)=Q$.

If $i<\infty$, choose any element $s\in\fp_i\setminus\fp_{i+1}$ so that we have $V(x_1R_{\fp_i})=\{\fp_i\}$ in $\Spec(R_{\fp_i})$.
Then the intersection $\bigcap_{n\ge1}(s^n)$ is a prime ideal by \cite[Lemma II.1.3]{fs01}. Notice that for every $n$ we must have  $\fp_{i+1}\subseteq (s^n)\subseteq (s)$. Otherwise, the total ordering of the ideals implies $s^n\in\fp_{i+1}$ and thus $s\in \fp_{i+1}$; a contradiction. Hence, $\fp_{i+1}= \bigcap_{n\ge1}(s^n)$, as there are no prime ideals in between $\fp_{i+1}\subset \fp_i$. Now consider $R_{\fp_{i+1}}=(R_{\fp_{i}})_{\fp_{i+1}}$, where we invert
\begin{equation*}
   S=R_{\fp_{i}}\setminus \fp_{i+1}R_{\fp_{i}}=R_{\fp_{i}}\setminus \bigcap\limits_{n\ge1}s^nR_{\fp_{i}}=\bigcup\limits_{n\ge1}(R_{\fp_{i}}\setminus s^nR_{\fp_{i}}). 
\end{equation*}
In particular, we may replace $S$ by $T=\bigcup_{n\ge2}(sR_{\fp_{i}}\setminus s^nR_{\fp_{i}})$, as $S^{-1}R_{\fp_{i}}=T^{-1}R_{\fp_{i}}$. Recall that $T^{-1}R_{\fp_{i}}\cong \colim_{t\in T}\, \frac{1}{t}R_{\fp_{i}}$, where the system runs over the directed system given by morphisms $\frac{1}{t}R_{\fp_{i}}\to\frac{1}{t'}R_{\fp_{i}}$ which exist if and only if there exists some $u\in T$ such that $t^\prime=tu$. By the description of $T$ we see that the powers of $s$ are final in this system, i.e., $R_{\fp_{i+1}}\cong S^{-1}R_{\fp_{i}}=R_{\fp_{i}}\big[\frac{1}{s}\big]$.
Finally, the conclusion follows from \Cref{prop:g(P)-in-D(R)} applied to $n=1, x_1=s$ since then
$g(\sfP_i)\simeq(R_{\fp_i}\to R_{\fp_{i+1}})$.
\end{proof}

\begin{proposition}\label{prop:BFreduction}
In the notation of \Cref{lem:BF} we have $\loc(g(\sfP_i))=\loc(k(\fp_i))$. 
\end{proposition}
\begin{proof}
The case $i=\infty$, i.e., $\loc(g(\sfP_i))=\loc(Q)$ follows from \Cref{lem:BF}.  For a fixed $i<\infty$ we also obtain $\loc(g(\sfP_i))=\loc(R_{\fp_{i+1}}/R_{\fp_{i}})$, as $g(\sfP_i)\simeq R_{\fp_{i+1}}/R_{\fp_{i}}[-1]$.

By \Cref{thm:propertiesStrong}(ii) we have $\fp_{i}R_{\fp_{i}}=sR_{\fp_{i}}$ for some non-zero $s\in R_{\fp_{i}}$. Without loss of generality, we may assume $s=s/1$ for some $s\in \fp_i\setminus\fp_{i+1}$. 
Then, by the same argument as in the proof of \Cref{lem:BF}, $\fp_{i+1}= \bigcap_{n\ge1}(s^n)$ and $R_{\fp_{i+1}}=R_{\fp_{i}}\big[\frac{1}{s}\big]$. Now the union of the filtration
\[
R_{\fp_i}\subset \frac{1}{s}R_{\fp_{i}}\subset \frac{1}{s^2}R_{\fp_{i}}\subset \cdots,
\]
whose terms are isomorphic to  $R_{\fp_{i}}$ is precisely $R_{\fp_{i}}\big[\frac{1}{s}\big]=R_{\fp_{i+1}}$. 
Passing the filtration to quotient modules, we get
\begin{equation}\label{eq:locFiltration}
\begin{aligned}
   \left(\frac{1}{s}R_{\fp_{i}}\right)/R_{\fp_{i}}\subset \left(\frac{1}{s^2}R_{\fp_{i}}\right)/R_{\fp_{i}}\subset \cdots\,\&\,\bigcup\limits_{n\ge1} \left[\left(\frac{1}{s^n}R_{\fp_{i}}\right)/R_{\fp_{i}}\right]=R_{\fp_{i+1}}/R_{\fp_{i}}.
   \end{aligned}
\end{equation}
This gives us a filtration of the module $R_{\fp_{i+1}}/R_{\fp_{i}}$ we are interested in. Now we wish to relate this to $k(\fp_i)$. To this end, notice that for any $n\ge1$ the composite $s^{-n}R_{\fp_{i}}\ontop{$\cdot s^n$}{\to} R_{\fp_{i}}\to R_{\fp_{i}}/s^nR_{\fp_{i}}$ is surjective and has kernel $R_{\fp_{i}}$. In particular, $s^{-n}R_{\fp_{i}}/R_{\fp_{i}}\cong R_{\fp_{i}}/s^nR_{\fp_{i}}$ as $R_{\fp_{i}}$ modules, which implies they are also isomorphic as $R$-modules. A similar argument reveals that multiplication by $s^{1-n}$ induces an isomorphism $s^{n-1}R_{\fp_{i}}/s^nR_{\fp_{i}}\cong R_{\fp_{i}}/ sR_{\fp_{i}}=k(\fp_i)$ of $R$-modules. Jointly, these modules participate in a short exact sequence
\[
\begin{tikzcd}[row sep=0.5cm,column sep=0.5cm]
0 \arrow[r] & s^{n-1}R_{\fp_{i}}/s^nR_{\fp_{i}} \arrow[r] & R_{\fp_{i}}/s^nR_{\fp_{i}} \arrow[r, "\varphi"] & R_{\fp_{i}}/s^{n-1}R_{\fp_{i}} \arrow[r] & 0,
\end{tikzcd}
\]
which witnesses $R_{\fp_{i}}/s^nR_{\fp_{i}}\cong s^{-n}R_{\fp_{i}}/R_{\fp_{i}}$ as an $n-$fold extension of $k(\fp_i)\cong\mathrm{ker}(\varphi)$.
Similarly, multiplication by $s$ induces an isomorphism $R_{\fp_{i+1}}/R_{\fp_{i}}\cong R_{\fp_{i+1}}/sR_{\fp_{i}}$ and these modules form an exact sequence
\[
\begin{tikzcd}[row sep=0.5cm,column sep=0.5cm]
0 \arrow[r] & R_{\fp_{i}}/sR_{\fp_{i}} \arrow[r] & R_{\fp_{i+1}}/sR_{\fp_{i}} \arrow[r, "\psi"] & R_{\fp_{i+1}}/R_{\fp_{i}} \arrow[r] & 0.
\end{tikzcd}
\]
If follows that $k(\fp_i)$ is the kernel of a surjective endomorphism of $R_{\fp_{i+1}}/R_{\fp_{i}}$.
Finally, since localising subcategories are closed under extensions, cokernels of monomorphisms of $R$-modules, and countable direct limits, we use \eqref{eq:locFiltration} and the exact sequences to conclude 
\[
\loc(R_{\fp_{i+1}}/R_{\fp_{i}})\subseteq 
\loc\left(\left(\frac{1}{s^n}R_{\fp_{i}}\right)/R_{\fp_{i}}\mid n\ge1\right)\subseteq \loc(k(\fp_i))\subseteq\loc(R_{\fp_{i+1}}/R_{\fp_{i}}). 
\]
But since $\loc(g(\sfP_i))=\loc(R_{\fp_{i+1}}/R_{\fp_{i}})$, the desired result follows.
\end{proof}
\begin{proposition}\label{thm:strongStratified}
Let $R$ be strongly discrete. Then $\sfD(R)$ is stratified by $\Spc(\sfD(R)^\omega)$.
\end{proposition}
\begin{proof}
The minimality at $\sfP_i$ holds because $\loc(k(\fp_i))$ is minimal and equals $\loc(g(\sfP_i))$ by \Cref{prop:BFreduction}. Then by \Cref{prop:l2gForStrongDiscrete} the stratification follows.
\end{proof}
\begin{theorem}\label{cor:classificationOfStratificationForVD}
Let $R$ be a valuation domain. Then the following are equivalent:
\begin{enumerate}
    \item $R$ is a strongly discrete valuation domain.
    \item $\sfD(R)$ is stratified by $\Spc(\sfD(R)^\omega)$.
    \item $\sfD(R)$ satisfies the Telescope Conjecture.
\end{enumerate}
\end{theorem}
\begin{proof}
By \cite[\S7]{gldm1} the Telescope Conjecture holds if and only if $R$ is strongly discrete. Then \Cref{thm:strongStratified} yields stratification.
If, on the other hand, the stratification holds, then $\Spec(R)$ is weakly noetherian by definition, and consequently noetherian by~\Cref{lem:T_D-for-VD}.
So the stratification implies the Telescope Conjecture by \cite[Theorem~9.11]{bhs23}.
\end{proof}
\begin{remark}
\Cref{cor:classificationOfStratificationForVD} gives a new proof that the Telescope Conjecture holds for $\sfD(R)$, where $R$ is a strongly discrete valuation domain, independent of \cite{gldm1}.
Indeed, \Cref{thm:strongStratified} shows that $\sfD(R)$ is stratified and, since $\Spec(R)\simeq\Spc(\sfD(R)^\omega)$ is a noetherian topological space, the Telescope Conjecture also follows from~\cite[Theorem 9.11]{bhs23}.
In forthcoming work \cite{BHSZ}, Barthel--Heard-Sanders-Zou provide an  independent proof of \Cref{cor:classificationOfStratificationForVD}.
\end{remark}

\section{The Gabriel dimension}
Instead of arguing by elementary means as in the previous section, we shall provide an alternative proof of \Cref{thm:strongStratified} using the Gabriel dimension. We will briefly recall the concept here and refer to \cite[\S5.5]{pop73}, \cite{gab62} for more background.

\begin{definition}\label{def:Gabriel-dim}
Let $\cC$ be a Grothendieck category. We define a sequence of localising (i.e.\ coproduct closed Serre) subcategories $\mathbb{L}_\alpha\subseteq\cC$ indexed by ordinal numbers as follows:
\begin{itemize}
\item We let $\mathbb{L}_0\subseteq\cC$ be the class of all semiartinian objects.
\item If $\alpha+1$ is a successor ordinal, $\mathbb{L}_{\alpha+1}$ is defined as the class of all objects of $\cC$ that are semiartinian in the localisation $\cC/\mathbb{L}_\alpha$.
\item For a limit ordinal $\beta$, we let $\mathbb{L}_\beta$ be the smallest localising subcategory containing $\bigcup_{\alpha<\beta}\mathbb{L}_\alpha$.
\end{itemize}
The \emph{Gabriel dimension} of $\cC$ is defined to be the least ordinal number $\beta$ such that $\cC=\mathbb{L}_\beta$, and is undefined if no such ordinal exists.
\end{definition}

If $\cC$ is a locally coherent Grothendieck category, it is natural to ask how the filtration of the Gabriel dimension translates to filtrations of Serre subcategories in the full subcategory $\cC^{\mathrm{fp}}$ formed by the finitely presented objects.
This is discussed in  \cite{bur01}, where Burke introduces the so-called f.p.-dimension.
\begin{definition}\label{def:fp-dim}
If $\cC$ is a locally coherent Grothendieck category, we can construct a similar filtration for the abelian category $\cC^{\mathrm{fp}}$ of finitely presentable objects.
\begin{itemize}
\item We let $\mathbb{F}_0\subseteq\cC^{\mathrm{fp}}$ be the class of objects of finite length.
\item If $\alpha+1$ is a successor ordinal, $\mathbb{F}_{\alpha+1}$ is defined as the class of all objects of $\cC^{\mathrm{fp}}$ that are of finite length in the Serre quotient $\cC^{\mathrm{fp}}/\mathbb{F}_\alpha$.
\item For a limit ordinal $\beta$, we let $\mathbb{F}_\beta=\bigcup_{\alpha<\beta}\mathbb{F}_\alpha$.
\end{itemize}
The \emph{f.p.-dimension} of $\cC$ is defined to be the least ordinal number $\beta$ such that $\cC^{\mathrm{fp}}=\mathbb{F}_\beta$, and is undefined if no such ordinal exists.
\end{definition}

\begin{remark}\label{rem:fp-dim}
There is a classical one-to-one correspondence between localising subcategories of finite type of $\cC$ and Serre subcategories of $\cC^{\mathrm{fp}}$; see~\cite[Theorem~2.8]{Herzog:1997}
or~\cite[Corollary~2.10]{loc_coh_grothendieck_spectrum}.
In particular, each Serre subcategory $\mathbb{F}_\alpha\subseteq\cC^{\mathrm{fp}}$ extends by direct limits to a localising sucategory $\mathbb{H}_\alpha\subseteq\cC$ (keeping the notation of~\cite{bur01}).
So $\cC$ has f.p.-dimension $\alpha$ if and only if $\alpha$ is the smallest ordinal such that $\cC=\mathbb{H}_\alpha$.
\end{remark}

A closely related concept is also that of $m$-dimension of a lattice and of a locally coherent Grothendieck category.

\begin{definition}[{\cite[\S7.2]{pp}}]\label{def:m-dim}
Given a modular lattice $L$ with the top and the bottom element, we denote by $\sim_L$ the smallest congruence on $L$ such that $a\sim_L b$ whenever $a<b$ and there is no element $c\in L$ strictly between $a$ and $b$.

By transfinite induction, we can now define congruences $\sim^\alpha$ and lattice quotients $L_\alpha=L/\sim^\alpha$ as follows:
\begin{itemize}
\item $L_0\coloneqq L$ and $\sim^0$ is the trivial congruence.
\item If $\alpha+1$ is a successor ordinal, we put $L_{\alpha+1} = L_\alpha/\sim_{L_\alpha}$ and $\sim_{\alpha+1}$ is the congruence corresponding to the composite quotient map
\[ L\lra {L/\sim_\alpha}=L_\alpha\lra L_{\alpha+1}={L_\alpha/\sim_{L_\alpha}}. \]
\item If $\beta$ is a limit ordinal, we set $\sim_\beta=\bigcup_{\alpha<\beta}\sim_\alpha$.
\end{itemize}
The \emph{m-dimension} of $L$ is defined to be the least ordinal number $\beta$ such that $L_{\beta+1}$ is the trivial (one-point) lattice, or undefined if no such ordinal exists. If $L=\{*\}$, we set it to $-1$.
\end{definition}

\begin{remark}\label{rem:m-dim}
If $L\neq\{*\}$, the definition implicitly uses the fact that the smallest ordinal $\gamma$ such that $L_\gamma=\{*\}$ is always a successor ordinal; see~\cite[\S7.1, p.~282]{pp}.
\end{remark}

Now we can summarise a key result of \cite{bur01} as follows.
\begin{theorem}[{\cite[Theorems 4.5 \& 5.1]{bur01}}]\label{thm:kgforcoherent}
    Let $\cC$ be a locally coherent Grothendieck category. Then the following are equivalent:
    \begin{enumerate}
        \item $\cC$ has Gabriel dimension.
        \item $\cC$ has f.p.-dimension.
        \item $\cC$ has m-dimension, i.e., for every $X\in\cC^{\mathrm{fp}}$ the modular lattice of finitely presented subobjects, denoted by $\lfp(X)$, has m-dimension.
    \end{enumerate}
\end{theorem}

Let $R$ be a coherent ring, e.g., a valuation domain, and recall that $\cC=\MOD_R$ is a locally coherent Grothendieck category. Then we have:

\begin{corollary}\label{prop:kgforcoherentrings}
    Let $R$ be coherent. Then the following are equivalent:
    \begin{enumerate}
        \item $\MOD_R$ has Gabriel dimension.
        \item $\lfp(R)$ has m-dimension.
        \item There exists no embedding of $(\mathbb{Q},\le)\hookrightarrow\lfp(R)$, i.e., there are no dense chains in $\lfp(R)$.
    \end{enumerate}
\end{corollary}
\begin{proof}
    The connection between m-dimension and dense chains is standard and can be found in \cite[Proposition 7.2.3]{pp} or \cite[Lemma 14.1.11]{kra22}. Now assume $\lfp(R)$ has m-dimension. The case when this dimension is $-1$ is clear. Otherwise, $R$ must have non-zero finite length in $\MOD_R/\mathbb{H}_{\alpha}$ for some ordinal $\alpha$ by \cite[Theorem 4.5]{bur01}. In particular, by definition of the filtration $R\in\mathbb{H}_{\alpha+1}$, and thus $\MOD_R$ has f.p.-dimension, hence, Gabriel dimension. The converse is \Cref{thm:kgforcoherent}.
\end{proof}

Getting back to our alternative approach to \Cref{thm:strongStratified}, we are going to use a key result of \cite{isolated}.

\begin{theorem}[{\cite[Theorem 2.1]{isolated}}, {\cite[Theorem 9.6]{bhsz26}}]\label{thm:kgimpliesstratification}
    Let $R$ be a commutative ring such that $\Spec(R)^\vee$ is a $T_D$-space. Then, if $\MOD_R$ has Gabriel dimension, $\sfD(R)$ is stratified by $\Spc(\sfD(R)^\omega)$.
\end{theorem}
We are now in the position to show that for valuation domains $R$, Gabriel dimension of $\MOD_R$ exists precisely when $R$ is strongly discrete. Hence, we provide a converse to \Cref{thm:kgimpliesstratification} in case $R$ is a valuation domain.

\begin{corollary}\label{cor:kgforvd}
    Let $R$ be a valuation domain. Then $\MOD_R$ has Gabriel dimension if and only if $R$ is strongly discrete.
\end{corollary}
\begin{proof}
    Using \Cref{prop:kgforcoherentrings}, this reduces to the question of the existence of dense chains of the lattice of (non-zero) finitely generated ideals. Such ideals are principal, and so are classified by non-zero elements of $R$ up to multiplication by a unit of $R$, and so also by elements of the positive cone $\Gamma(R)^+=\{v\in\Gamma(R)\mid v\ge 0\}$ of the value group $\Gamma(R)$ of $R$ (cp.~\cite[Proposition 2.8]{bt25}). Thus, it suffices to show: $R$ is strongly discrete if and only if $(\mathbb{Q},\leq)$ does not embed into $(\Gamma(R),\leq)$. By  \cite[Corollary~4.2]{puninskii_superdecomposable} the latter is equivalent to the non-existence of superdecomposable pure-injective $R$-modules. Then \cite[Theorem XIII.5.11]{fs01} yields the result.
\end{proof}
\begin{theorem}\label{thm:kg=stratforvd}
    Let $R$ be a valuation domain. Then the following are equivalent:
\begin{tasks}[style=enumerate, label-align=right ,column-sep = {-7mm},item-indent = {7mm},label-offset = {0.8em} ](2)
        \task*(1)[(i)] $R$ is strongly discrete.
        \task*[(ii)] $\sfD(R)$ satisfies the Telescope Conjecture.
        \task*(1)[(iii)] $\sfD(R)$ is stratified by $\Spc(\sfD(R)^\omega)$.
        \task*[(iv)] $\MOD_R$ has Gabriel dimension.
    \end{tasks}
\end{theorem}
\begin{proof}
   Again, $(i)$ and $(ii)$ are equivalent by \cite[\S7]{gldm1}, and $(i)$ and $(iv)$ are equivalent by \Cref{cor:kgforvd}.
   If $(i)$ holds, then $\Spec(R)^\vee$ is a $T_D$-space by \Cref{thm:propertiesStrong} and \Cref{lem:T_D-for-VD}, so it follows from \Cref{thm:kgimpliesstratification} that $(i)$ implies $(iii)$.
   Finally, we can use the same argument as in~\Cref{cor:classificationOfStratificationForVD} for the fact that $(iii)$ implies $(ii)$. If $(iii)$ holds, then $\Spec(R)$ is weakly noetherian, so noetherian by \Cref{lem:T_D-for-VD}, and then $(ii)$ follows from \cite[Theorem~9.11]{bhs23}.
\end{proof}
In particular, by \Cref{thm:kg=stratforvd} there is an isomorphism $\Loc(\sfD(R))\cong \mathcal{P}(\Spec(R))$ for strongly discrete valuation domains $R$ . However, the following example shows that the mere existence of such an isomorphism for a given valuation domain $R$ does not necessarily imply that $R$ is strongly discrete.
\begin{example}\label{ex:ce}
    By \cite[Proposition 5.1, Remark 5.2]{bt25} there exists a valuation domain $R$ with the following properties:
\begin{enumerate}
    \item The value group is the lexicographic product of copies of $\Z$, indexed by $\N$. 
    \item The prime ideals $\fp_i$ are indexed by $\N\cup\{\infty\}$, form an ascending chain by the natural order, $\fp_0=(0)$, and $\fp_\infty=\bigcup_{i\in\N}\fp_i$ is maximal.
    \item The only idempotent prime ideals are zero and $\fp_\infty$.
\end{enumerate}
    As $\fp_\infty$ is idempotent and non-zero, $\sfD(R)$ cannot satisfy the Telescope Conjecture by \cite[\S7]{gldm1}, hence Gabriel dimension for $\MOD_R$ is undefined by \Cref{thm:kg=stratforvd}. Moreover,
    $\Spec(R)\simeq\Spc(\sfD(R)^\omega)$ is not weakly noetherian by \Cref{lem:T_D-for-VD} (in fact $\fp_\infty\in\Spec(R)$ is not weakly visible).

    However, we claim that  $\sfL\coloneqq\loc(k(\fp_i)\mid i\in\N\cup\{\infty\})$ equals $\sfD(R)$, and therefore $\Loc(\sfD(R))\cong\mathcal{P}(\Spec(R))$. Thus, $R$ provides a counterexample to the converse of \cite[Theorem 2.1]{isolated} and shows that having Gabriel dimension is not necessary for residue fields to generate.
    
    To prove $\sfL=\sfD(R)$, let $i\in\N$ and define $\sfL_i\coloneqq\loc(k(\fp)\mid \fp\subseteq \fp_i)\subseteq\sfD(R)$. As $R_{\fp_i}$ has discrete value group and is of finite Krull dimension, it is strongly discrete \cite[Proposition~II.8.4]{fs01}. Hence, using \Cref{prop:l2gForStrongDiscrete} and \Cref{prop:BFreduction}, (or alternatively \Cref{thm:kg=stratforvd} together with \cite[Theorem 2.1]{isolated}), we find that $\sfD(R_{\fp_i})$ is generated by residue fields. As the restriction functor $\sfD(R_{\fp_i})\to \sfD(R)$ is fully faithful and coproduct-preserving, we obtain $\sfL_i=\loc(R_{\fp_i})$. Since $\sfL_i\subset \sfL$, we have $R_{\fp_i}\in \sfL$ for any $i\in\N$. Now, crucially, $\fp_i R_{\fp_i} = \fp_i$ as subsets of $R_{\fp_i}$ since $\fp_i$ is canonically an $R_{\fp_i}$-module by~\cite[Lemma~4.3(iv)]{fs01}, and $\fp_iR_{\fp_i}=s_iR_{\fp_i}$ for some $s_i\in R$ since $R_{\fp_i}$ is strongly discrete.
    Hence, $\fp_i\simeq R_{\fp_i}$ as $R$-modules for $i\neq0$ and so $\fp_i\in\sfL$. As $\fp_\infty$ is the countable filtered colimit of $\fp_i$, we have $\fp_\infty\in\sfL$. As $k(\fp_\infty)=R/\fp_\infty\in\sfL$, this forces $R\in\sfL$ and $\sfL=\sfD(R)$. This proves the claim.

    We have now proven that the small Foxby support stratifies in this case, which replaces  the Balmer-Favi support in \Cref{thm:stratification-abstract} and realizes the isomorphism $\Loc(\sfD(R))\cong\mathcal{P}(\Spec(R))$.
    Further, we will provide a negative answer saying that the smashing support as defined in~\cite[\S3.2]{bigspectra}, that is strictly finer than the Balmer-Favi support in our example, does not do the job for stratification.

    In order to analyse the smashing support, we first need to compute the lattice $\Sm(\sfD(R))$ of smashing localising subcategories of $\sfD(R)$, which is known to be a frame (i.e.\ a complete lattice where finite meets distribute over arbitrary joins) by the main result of~\cite{bks20}.
    The following construction is due to \cite[\S5]{gldm1}, but for ease we will refer to \cite[Construction 3.3, Theorem 3.4]{bt25}. For our ring $R$, let $j\ge-1$. The elements of the frame, which are representatives of equivalence classes of homological epimorphisms originating from $R$, can be described as follows:

    \begin{eqnarray}\label{eqn:epis}
    \begin{split}
              \mbox{Flat epis:}\enspace &f_{-1}\colon R\to 0,\enspace f_\infty\colon R\to R, \enspace f_j\colon R\to R_{\fp_j} &\mbox{ for }j\ge0.\\
        \mbox{Non-flat epis:}\enspace&g_{-1}\colon R\to k(\fp_\infty),\enspace\enspace g_j\colon R\to R_{\fp_j}\dert_R k(\fp_\infty) &\mbox{ for }j\ge0.  
    \end{split}
    \end{eqnarray}
    A homological epimorphism $f\colon R\to S$ corresponds to the smashing localising subcategory $\sfS=\ker(S\dert_R-)\subseteq\sfD(R)$. Adopting notation of \cite[\S 4.1]{bt25},
    we obtain \Cref{fig:smash}.
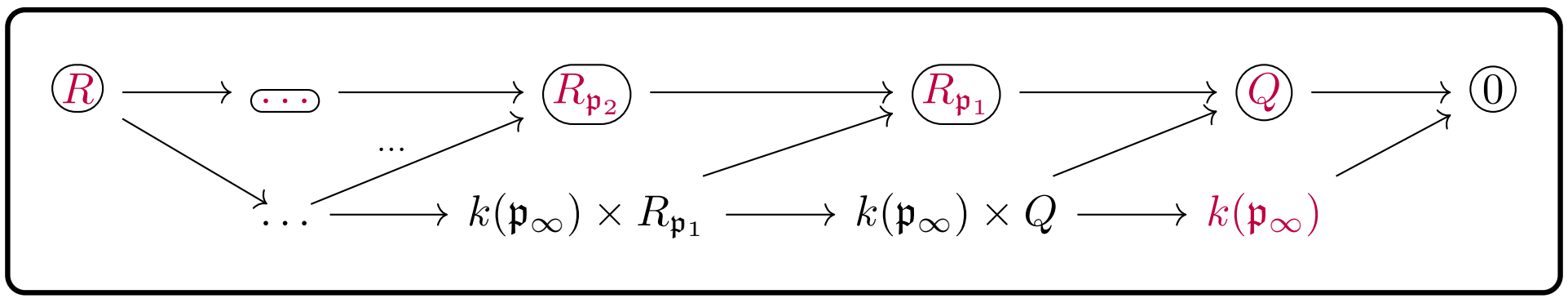
\begin{figure}[h]
\begin{center}
\begin{tikzcd}[style={every outer matrix/.append style={draw=black, inner xsep=6pt , inner ysep=9pt, rounded corners, very thick}},row sep=0.3cm,column sep=0.9cm]

\Circled[inner color=purple]{R} \arrow[r] \arrow[rd]                      & \Circled[inner color=purple]{\hdots} \arrow[r]                    & \Circled[inner color=purple]{R_{\fp_2}} \arrow[r]                                & \Circled[inner color=purple]{R_{\fp_1}} \arrow[r]                        & \Circled[inner color=purple]{Q} \arrow[r]              & \Circled{0} \\
 & \hdots \arrow[r] \arrow[ru, "\hdots" ] & k(\fp_\infty)\times R_{\fp_1} \arrow[r] \arrow[ru] & k(\fp_\infty)\times Q \arrow[r] \arrow[ru] & {\color{purple}k(\fp_\infty)} \arrow[ru] &  
\end{tikzcd}
\begin{center}\caption{The Hasse diagram of the smashing frame $\Sm(\sfD(R))$ in terms of the codomains of the representing homological epimorphisms \Cref{eqn:epis}. Arrows correspond to natural factorisations, which also correspond to inclusions in $\Sm(\sfD(R))$. The compactly generated ones, i.e., flat ones, are encircled, while meet-primes in the frame are coloured.}\label{fig:smash}\end{center}
\end{center}
\end{figure} 
Note that the bottom element $0=\ker(R\dert_R-)$ of the smashing frame is a meet-prime as the intersection of two non-zero smashing subcategories is non-zero. The compactly generated localisations of $\sfD(R)$ correspond to the flat epimorphisms and are encircled in the figure. Let $\sfC_j\in\Sm(\sfD(R))$ correspond to $f_j$ for $-1\le j\le \infty$, which is the upper row in \Cref{fig:smash}, and analogously $\sfS_j$ to $g_j$ for $-1\le j<\infty$; the lower row.
By the kernel description and $\Loc(\sfD(R))\cong\cP(\Spec(R))$ we get
\[\
\sfS_j=\loc(k(\fp_i)\mid j<i<\infty),\,\sfC_\infty=0,\, \sfC_j=\loc(k(\fp_i)\mid j<i\le \infty)\mbox{ for }-1\le j<\infty.
\]

Having described the smashing frame, we can construct the smashing spectrum $\Spc^s(\sfD(R))$ that we obtain by Stone duality; we refer to \cite[\S2.1 and \S3]{bigspectra} for details.
The points of $\Spc^s(\sfD(R))$ are the meet-primes of $\Sm(\sfD(R))$, i.e.\ the coloured points in \Cref{fig:smash}.
We also have the smashing support: if $x\in\sfD(R)$, then $\supp^s(x)=\{\sfP\in\Spc^s(\sfD(R)) \mid x\notin\sfP\}$. An open basis of $\Spc^s(\sfD(R))$ is given by the sets
\[
U_{\sfS} = \{ \sfP\in\Spc^s(\sfD(R)) \mid \sfS\not\subseteq\sfP \},
\]
where $\sfS$ runs over the elements of the smashing frame, i.e.\ the smashing localising subcategories of $\sfD(R)$. Note also that for each $\sfS\in\Sm(\sfD(R))$, the following holds tautologically by~\cite[Remark~3.2.5]{bigspectra}:
\[ \supp^s(\sfS) \coloneqq \bigcup_{x\in\sfS}\supp^s(x) = U_{\sfS}. \]

We claim that $\supp^s$ cannot induce an isomorphism of lattices given by mapping a localising subcategory to its support,
\begin{align*}
\Loc(\sfD(R)) &\simeq \cP(\Spc^s(\sfD(R))), \\
\sfL &\mapsto \supp^s(\sfL)\coloneqq\bigcup_{x\in\sfL}\supp(x).
\end{align*}
To this end, notice that by our explicit computations, we have
\[
\supp^s(\sfS_j)=\{\sfP\in\Spc^s(\sfD(R))\mid\sfS_j\not\subseteq\sfP\}=\{\sfC_\infty,\sfC_\ell\mid j<\ell\}.
\]
Furthermore, $\bigcap_{j\ge0}\sfS_j=\sfC_\infty$, while $\bigcap_{j\ge0}\supp^s(\sfS_j)=\{\sfC_\infty\}$. Hence, if $\supp^s$ induced an order isomorphism, and thus stratified $\sfD(R)$ by $\Spc^s(\sfD(R))$, it would preserve these meets. This is absurd, however, as we would obtain
\[
\varnothing=\supp^s(\sfC_\infty)=\supp^s(\bigcap_{j\ge0}\sfS_j)=\bigcap_{j\ge0}\supp^s(\sfS_j)=\{\sfC_\infty\}.
\]
\end{example}

\section{Stratification for rings of weak global dimension at most one}
In this section we characterise for a certain class of commutative rings of weak global dimension at most one when $\Spc(\sfD(R)^\omega)$ stratifies $\sfD(R)$. 
Recall from \Cref{lem:wgld1} that a commutative ring has weak global dimension at most one if and only if the localisation at any prime is a valuation domain.
By \cite[\S 7]{gldm1} we know that for such rings, $\sfD(R)$ satisfies the Telescope Conjecture precisely when $R_\fp$ is a strongly discrete valuation domains for all prime ideals $\fp$.

\begin{definition}[Cf. {\cite[Proposition III.7.4]{fs01}}]\label{def:strongp}
    A commutative ring $R$ of weak global dimension at most one is called {\it strongly discrete} if for any prime $\fp$ the localisation $R_\fp$ is a strongly discrete valuation domain. 
\end{definition}
Recall from \Cref{thm:propertiesStrong} that the ascending chain condition on primes holds for strongly discrete valuation domains.
The good news is:
\begin{lemma}\label{lem:weaklynoethprüfer}
Let $R$ be strongly discrete. Then $\Spc(\sfD(R)^\omega)^\vee$ is a $T_D$-space.
\end{lemma}
\begin{proof}
    We shall prove that for any prime ideal $\fp$ there are Thomason subsets $U,V$ in $\Spec(R)$ satisfying $U\cap V^c=\{\fp\}$.
    As $R$ is strongly discrete, there exists by \Cref{lem:wgld1} and  \Cref{thm:propertiesStrong} some $x\in \fp$ such that $x R_\fp=\fp R_\fp$, as $R_\fp$ is local. Note that $\fp$ is minimal over $(x)$. Let $V\coloneqq V(x)\setminus\{\fp\}$ and $\fq\in V$. Then by construction, we cannot have $\fq\subsetneq \fp$. Thus, there exists some $y\in\fq\setminus\fp$, and we have $(x,y)\not\subseteq \fp, V(x,y)\subseteq V$. Therefore, we may choose the Thomason subsets $\bigcup_{\fq\in V,y\in\fq\setminus\fp}V(x,y)=V$ and $U=V(x)$, which satisfy $U\cap V^c=\{\fp\}$.
%
\end{proof}
Now, the bad news is that there is no analogous statement of \Cref{cor:classificationOfStratificationForVD} for strongly discrete rings of weak global dimension at most one.
This is because it already fails for absolutely flat rings, i.e. rings of weak global dimension zero. Explicitly for a field $K$, the infinite product $\prod_\N K$ serves a counterexample; see \cite[Example 4.13]{stevenson_absolutelyflatrings}. However, we can rely on the following.
\begin{theorem}[{\cite[Theorem 7.7]{isolated}}]\label{thm:scattered}
    For a commutative ring $R$, the local-to-global principle holds for $\sfD(R)$ if and only if $\Spec(R)^\vee$ is scattered.
\end{theorem}
\begin{remark}\label{rem:onlycusps}
    Recall that by \cite[Lemma 7.3]{isolated} $\Spec(R)^\vee$ is \emph{scattered} if and only if every non-empty closed subset $C$ of $\Spec(R)^\vee$ has an isolated point $\fp$, i.e., $\fp\not\in\overline{(C\setminus\{\fp\})}^{\vee}$. By \cite[Corollary 4.1.6]{dst19} any such $C$ is generalisation closed in $\Spec(R)$. In our case, in view of \Cref{def:strongp}, the only points that could potentially be isolated in $C$ are the maximal primes in $C$, i.e., those without a successor in $C$. 
\end{remark}
\begin{notation}
    Let $R$ be strongly discrete of weak global dimension at most one. For a non-empty generalisation-closed subset $C\subseteq\Spec(R)$ let $\Max(C)$ denote the maximal elements of $C$. Note that this is non-empty as $R$ is strongly discrete.
\end{notation}
\begin{lemma}\label{rem:sharp}
    Let $R$ be a strongly discrete ring of weak global dimension at most one and let $C$ be a non-empty closed subset of $\Spec(R)^\vee$. Then $\fp\in C$ is an isolated point in $C$ if and only if there exists a finitely generated ideal $I$ such that $I\subseteq \fp$, but $I\not\subseteq\fq$ for $\fp\neq\fq\in\Max(C)$. 
\end{lemma}
\begin{proof}
    Assume $\fp$ is isolated in $C$. Then \Cref{rem:onlycusps} gives $\fp\in\Max(C)$. By definition, there exists an open $U\subseteq\Spec(R)^\vee$ such that $C\cap U=\{\fp\}$. Since the $V(I)$, with $I$ finitely generated, form a basis, we may replace $U$ by some $V(I)$. Hence, $C\cap V(I)=\{\fp\}$, i.e., $\fp$ contains $I$ and no other $\fq\in\Max(C)$ contains $I$.
    
    Conversely, assume that $I$ is a finitely generated ideal such that $I\subseteq \fp$, but $I\not\subseteq\fq$ for $\fp\neq\fq\in\Max(C)$. 
    Let $t\in \fp$ such that $tR_\fp=\fp R_\fp$, which exists as $R$ is strongly discrete, and consider the finitely generated ideal $I+tR$. From the assumptions $C\setminus \{\fp\}=C\cap (\Spec(R)\setminus V(I+tR))$ is an intersection of two closed sets in $\Spec(R)^\vee$, so is itself closed in $\Spec(R)^\vee$.
    But then $\fp\not\in\overline{(C\setminus\{\fp\})}^{\vee}$, i.e., $\fp$ is isolated in $C$.
\end{proof}
The property of \Cref{rem:sharp} has been studied in the context of Prüfer domains.
\begin{definition}[{Cf. \cite[Theorem 4.1.6]{sharp}}]\label{def:sharp}
    Let $R$ be strongly discrete of weak global dimension at most one and let $\varnothing\neq C\subseteq\Spec(R)^\vee$ be closed. Call $\fp\in\Max(C)$ {\it sharp in }$C$ if there exists a finitely generated ideal $I$ such that $\fp\in V(I)$ and $\fq\not\in V(I)$ for all $\fq\in\Max(C)$ with $\fq\neq\fp$.
\end{definition}
Let us provide a non-trivial example where the local-to-global principle holds.
\begin{example}
    This ring is due to \cite[\S 8.2]{sharp}. For a prime number $p$ let $V=\Z_{(p)}$. Then one can find a sequence of algebraic numbers $t_i$ such that $\Q_i\coloneqq\Q(t_1,\hdots,t_i)$ is a tower of field extensions and $t_{i+1}$ is quadratic over $\Q_i$. Then the integral closure $R$ of $V$ in $\Q_\infty$ is an almost Dedekind domain, which is not Dedekind, see Theorem 8.2.2 loc.\ cit. In other words, for any maximal ideal $\fm\subset R$, the localisation $R_\fm$ is a discrete rank one valuation domain, but $\Spec(R)$ is not noetherian. In particular, $R$ is a strongly discrete Prüfer domain. In fact, the proof of loc. cit. reveals that any maximal ideal $\fm$ is principal, except one. Call this single maximal ideal $\fm_0$. By \cite[Proposition 2.2]{noethSpec} $\Spec(R)$ is not noetherian precisely because of $\fm_0$, as it is not the radical of a finitely generated ideal. As all other $\fm\in\maxspec(R)$ are principal, we have $V(\fm)=\{\fm\}$, i.e., they are sharp in any non-empty closed $C\subseteq\Spec(R)^\vee$ containing them. On the other hand, $\fm_0$ is not sharp in $\Spec(R)^\vee$, as for any finitely generated ideal $I$ we have $V(I)\neq V(\fm_0)$.
    
    However, $\Spec(R)^\vee$ is still scattered.
    Indeed, if $C$ is non-empty closed and contains $\fm\in\maxspec(R)$ different from $\fm_0$, then $\fm$ is sharp since it is a principal maximal ideal.    
    The only other non-empty closed sets in $\Spec(R)^\vee$ are $C_1=\{(0)\}$ and $C_2=\{(0),\fm_0\}$. The case of $C_1$ is trivial, while $\fm_0$ is sharp in $C_2$ as $\Max(C_2)=\{\fm_0\}$.
\end{example}
Let us also discuss how sharpness can be detected by flat ring epimorphisms. By \cite[Proposition 6.7]{gldm1} there is a bijection between non-empty closed subsets in $\Spec(R)^\vee$ and equivalence classes of flat ring epimorphisms $f\colon R\to S\neq0$. The correspondence sends a representative $f$ to $C_S\coloneqq\{\fq\in\Spec(R)\mid \fq S\neq S\}$.
\begin{lemma}\label{lem:sharpepi}
    Let $R$ be strongly discrete of weak global dimension at most one and let $f\colon R\to S\neq0$ be a flat epimorphism. Then $\fp\in\Max(C_S)$ is sharp if and only if there exists $\fn\in\maxspec(S)$ and a finitely generated ideal $J\subseteq S$ such that $(S/J,\fn/J)$ is a local ring.
\end{lemma}
\begin{proof}
    As $f$ is a flat epimorphism of rings, $\Spec(f)$ is by~\cite[Corollaire~IV.2.2]{Lazard:1969} a homeomorphism onto its image, which is $C_S$. As such, it is straightforward to show that the condition of \Cref{def:sharp} translates to the existence of a finitely generated ideal $J\subseteq S$ such that $S/J$ is local.
    In particular, denoting by $\fn$ the unique maximal ideal of $V(J)\subseteq\Spec(S)$, $(S/J,\fn/J)$ is as required.

    Suppose conversely that $J\subseteq S$ is a finitely generated ideal and $(S/J,\fn/J)$ is a local ring. Since $S_\fn$ is a strongly discrete valuation domain, there exists $t\in S$ such that $\fn S_\fn=t S_\fn$. In particular, $V(J+tS)=\{\fn\}$ in $\Spec(S)\cong C_S$.
\end{proof}
\begin{definition}\label{def:sharpepi}
Let $R$ be a ring of weak global dimension at most one. We shall call the flat epimorphisms satisfying the condition in \Cref{lem:sharpepi} \emph{sharp}.
\end{definition}
Let us summarise our findings of this section in the following result.
\begin{theorem}\label{thm:wgld1strat}
    Let $R$ be a ring of weak global dimension at most one. Then the following are equivalent:
    \begin{enumerate}
        \item $R$ is strongly discrete and each non-zero flat epimorphism is sharp.
        \item $R$ is strongly discrete and the local-to-global principle holds.
        \item $\sfD(R)$ is stratified by $\Spc(\sfD(R)^\omega)$.
    \end{enumerate}
\end{theorem}
 \begin{proof}
    Suppose $R$ is strongly discrete and the local-to-global principle holds. Then $\Spc(\sfD(R)^\omega)^\vee$ is a $T_D$-space by \Cref{lem:weaklynoethprüfer} and $\sfD(R)$ is stratified at each local category $\sfD(R_\fp)$ by definition and \Cref{cor:classificationOfStratificationForVD}. Thus, $\sfD(R)$ is stratified by \cite[Corollary 5.3]{bhs23}. The other implications are straightforward given \Cref{thm:scattered} and \Cref{lem:sharpepi}.
 \end{proof}
Finally, using \Cref{thm:wgld1strat} we may extend the result of \Cref{thm:kg=stratforvd} to coherent rings of weak global dimension at most one, e.g., Prüfer domains.
\begin{corollary}\label{cor:kgandwgld}
    Let $R$ be a coherent ring of weak global dimension at most one, such that each non-zero flat epimorphism is sharp. Then the following are equivalent:
\begin{tasks}[style=enumerate, label-align=right ,column-sep = {-7mm},item-indent = {7mm},label-offset = {0.8em} ](2)
        \task*(1)[(i)] $R$ is strongly discrete.
        \task*[(ii)] $\sfD(R)$ satisfies the Telescope Conjecture.
        \task*(1)[(iii)] $\sfD(R)$ is stratified by $\Spc(\sfD(R)^\omega)$.
        \task*[(iv)] $\MOD_R$ has Gabriel dimension.
    \end{tasks}
\end{corollary}
\begin{proof}
Assume $R$ is strongly discrete. Using \Cref{thm:wgld1strat} and \Cref{thm:scattered}, $X=\Spec(R)^\vee$ is scattered. As usual, denote $X^{0}=X$, $X^{\alpha+1}=X^\alpha\setminus\mathrm{Isol}(X^\alpha)$ for successor ordinals $\alpha$, where $\mathrm{Isol}(X^\alpha)\subseteq X^\alpha$ denotes the set of isolated points, and $X^\lambda=\bigcap_{\alpha<\lambda}X^\alpha$ for limit ordinals $\lambda$. This procedure gives $X^\gamma=\varnothing$ for some ordinal $\gamma$, so for $\varnothing \neq U\in\Omega(X)^\omega$ there is a least ordinal $\beta$ with $U\cap X^\beta=\varnothing$. Note that $\beta$ cannot be a limit ordinal. Otherwise $(U\setminus(U\cap X^\alpha))_{\alpha<\beta}$ yields an open cover of $U$, thus there is a finite subcover, which contradicts $\varnothing=U\cap X^\beta=\bigcap_{\alpha<\beta}U\cap X^\alpha$.
    Hence, $\alpha(U)\coloneqq\beta-1$ is the largest ordinal satisfying $U_X\coloneqq U\cap X^{\alpha(U)}\neq\varnothing$. In particular, $U_X$ is closed in~$U$ and it must by finite be quasi-compactness of $U$.

    Let us now prove that $(i)$ implies $(iv)$.
    If Gabriel dimension of $\MOD_R$ is undefined, there exists a dense chain $\mathcal{C}$ of finitely generated ideals by \Cref{prop:kgforcoherentrings}. Pick   $I\subsetneq J$ in $\mathcal{C}$.
    Since $J/I\neq0$ is finitely generated, we have $\varnothing\ne\Supp(J/I)=V(\Ann(J/I))$.
    The annihilator of $J/I=m_1R+\cdots+m_kR$ is finitely presented by coherence, as it is the kernel of $R\to (J/I)^k$, $r\mapsto(m_1r,\dots,m_kr)$. In particular, $V(\Ann(J/I))=V(I\,\colon J)$ is Zariski closed with quasi-compact complement in $\Spec(R)$. 
    Thus, $\varnothing \neq \Supp(J/I)\in\Omega(X)^\omega$ and we may use the previous paragraph for any pair $I\subsetneq J$ in $\mathcal{C}$. To show some $R_\fp$ is not strongly discrete for some prime $\fp$, it suffices by \Cref{cor:kgforvd} to find $S\subsetneq T$ such that $(J/I)_\fp\neq0$ for any $S\subseteq I\subsetneq J\subseteq T$, where $S,I,J,T\in\mathcal{C}$. For ideals $I\subsetneq J$ in $\mathcal{C}$ define $f(I,J)\coloneqq (\alpha(\Supp(J/I)),\# \Supp(J/I)_X)$. Let the codomain of $f$ be its image and endow it with the lexicographical order. Out of this image, pick a pair, $S\subsetneq T$, in the preimage of the minimal value under this order. Let $S\subseteq I\subsetneq J\subseteq T$ be in $\mathcal{C}$, and set $U=\Supp(T/S)$. Then $\Supp(J/I)\subseteq U$ as $J/I$ is a subquotient of $T/S$. Hence, $\alpha (\Supp(J/I))\le\alpha(U)$, and by construction this forces equality. Thus, $\Supp(J/I)_X\subseteq U_X$. By construction we must have $\# U_X\le\# \Supp(J/I)_X$, and finiteness gives $\Supp(J/I)_X=U_X$, i.e., $f(I,J)=f(S,T)$. By design, we have now for each $\fp\in U_X$ that $\fp\in\Supp(J/I)$, i.e., $(J/I)_\fp\neq0$.
    This finishes the proof of the implication.

    Assume $\MOD_R$ has Gabriel dimension. Then $R$ must be strongly discrete by \Cref{thm:kg=stratforvd}, so $(i)$ and $(iv)$ are equivalent.
    
    If $R$ is strongly discrete then, with the assumption and \Cref{thm:wgld1strat}, the stratification of $\sfD(R)$ follows, showing that $(i)$ implies $(iii)$.

    If $\sfD(R)$ is stratified, $\Spec(R)$ must be weakly noetherian by definition, and so must be $\Spec(R_\fp)$ for each $\fp$. By \Cref{lem:T_D-for-VD}, this implies that $\Spec(R_\fp)$ is noetherian for each $\fp$. In particular, $\Spec(R)$ is generically noetherian in the sense of~\cite[Definition~9.5]{bhs23}, so~\cite[Theorem~9.11]{bhs23} applies, showing that the stratification implies the Telescope Conjecture.
    In other words, $(iii)$ implies $(ii)$.
    
    The final implication from $(ii)$ to $(i)$ follows by \cite[\S 7]{gldm1}, since strongly discrete rings $R$ are characterised among those of weak global dimension at most one by $\sfD(R)$ satisfying the Telescope Conjecture.
\end{proof}
\bibliography{bibliography}\bibliographystyle{alpha}

@incollection{BHSZ,
    AUTHOR = {Barthel, Tobias and Heard, Drew and Sanders, Beren and Zou,
              Changhan},
      title={The tensor triangular geometry of valuation domains},
            year={2026},
      journal={arXiv e-prints},
      note = "Forthcoming work"
}

@incollection{cosupport,
      title={Cosupport in tensor triangular geometry}, 
      author={Barthel, Tobias and Castellana, Natalia and Heard, Drew and Sanders, Beren},
      year={2023},
      journal={arXiv e-prints},
      eprint={2303.13480},
      archivePrefix={arXiv},
      primaryClass={math.CT},
      url={https://arxiv.org/abs/2303.13480}, 
      note = "arXiv preprint: \url{https://arxiv.org/abs/2303.13480}"
}

@incollection{bigspectra,
      title={Big categories, big spectra}, 
      author={Balchin, S. and Stevenson, G.},
      year={2023},
      journal={arXiv e-prints},
      eprint={2109.11934},
      archivePrefix={arXiv},
      primaryClass={math.AC},
      url={https://arxiv.org/abs/2109.11934}, 
      note = "arXiv preprint: \url{https://arxiv.org/abs/2109.11934}"
}

@article {ric97,
    AUTHOR = {Rickard, Jeremy},
     TITLE = {Idempotent modules in the stable category},
   JOURNAL = {J. London Math. Soc. (2)},
  FJOURNAL = {Journal of the London Mathematical Society. Second Series},
    VOLUME = {56},
      YEAR = {1997},
    NUMBER = {1},
     PAGES = {149--170},
      ISSN = {0024-6107,1469-7750},
   MRCLASS = {20J05 (20C05)},
  MRNUMBER = {1462832},
MRREVIEWER = {Amnon\ Neeman},
       DOI = {10.1112/S0024610797005309},
       URL = {https://doi.org/10.1112/S0024610797005309},
}

@article {tensoractions,
    AUTHOR = {Stevenson, Greg},
     TITLE = {Support theory via actions of tensor triangulated categories},
   JOURNAL = {J. Reine Angew. Math.},
  FJOURNAL = {Journal f\"ur die Reine und Angewandte Mathematik. [Crelle's
              Journal]},
    VOLUME = {681},
      YEAR = {2013},
     PAGES = {219--254},
      ISSN = {0075-4102,1435-5345},
   MRCLASS = {18E30 (18D10)},
  MRNUMBER = {3181496},
MRREVIEWER = {Pawel\ Sosna},
       DOI = {10.1515/crelle-2012-0025},
       URL = {https://doi.org/10.1515/crelle-2012-0025},
}

@incollection {noethSpec,
    AUTHOR = {Gabelli, Stefania},
     TITLE = {Generalized {D}edekind domains},
 BOOKTITLE = {Multiplicative ideal theory in commutative algebra},
     PAGES = {189--206},
 PUBLISHER = {Springer, New York},
      YEAR = {2006},
      ISBN = {978-0-387-24600-0; 0-387-24600-2},
   MRCLASS = {13F05 (13A15)},
  MRNUMBER = {2265809},
MRREVIEWER = {William\ W.\ Smith},
       DOI = {10.1007/978-0-387-36717-0\_12},
       URL = {https://doi.org/10.1007/978-0-387-36717-0_12},
}

@book {sharp,
    AUTHOR = {Fontana, Marco and Huckaba, James A. and Papick, Ira J.},
     TITLE = {Pr\"ufer domains},
    SERIES = {Monographs and Textbooks in Pure and Applied Mathematics},
    VOLUME = {203},
 PUBLISHER = {Marcel Dekker, Inc., New York},
      YEAR = {1997},
     PAGES = {x+328},
      ISBN = {0-8247-9816-3},
   MRCLASS = {13F05 (13A15 13G05)},
  MRNUMBER = {1413297},
MRREVIEWER = {Evan\ G.\ Houston},
}

@book {kra22,
    AUTHOR = {Krause, Henning},
     TITLE = {Homological theory of representations},
    SERIES = {Cambridge Studies in Advanced Mathematics},
    VOLUME = {195},
 PUBLISHER = {Cambridge University Press, Cambridge},
      YEAR = {2022},
     PAGES = {xxxiv+482},
      ISBN = {978-1-108-83889-4},
   MRCLASS = {16G10 (16E35 16E45 16G20 18E10 18G80 20C15)},
  MRNUMBER = {4327095},
MRREVIEWER = {E.\ R.\ Alvares},
}

@article {bhsz26,
    AUTHOR = {Barthel, Tobias and Heard, Drew and Sanders, Beren and Zou,
              Changhan},
     TITLE = {Homological stratification and descent},
   JOURNAL = {J. Inst. Math. Jussieu},
  FJOURNAL = {Journal of the Institute of Mathematics of Jussieu. JIMJ.
              Journal de l'Institut de Math\'ematiques de Jussieu},
    VOLUME = {25},
      YEAR = {2026},
    NUMBER = {2},
     PAGES = {1081--1125},
      ISSN = {1474-7480,1475-3030},
   MRCLASS = {18G80 (18F99 55P91 55U35)},
  MRNUMBER = {5023884},
       DOI = {10.1017/S1474748025101515},
       URL = {https://doi.org/10.1017/S1474748025101515},
}

@article {bur01,
    AUTHOR = {Burke, Kevin},
     TITLE = {Co-existence of {K}rull filtrations},
   JOURNAL = {J. Pure Appl. Algebra},
  FJOURNAL = {Journal of Pure and Applied Algebra},
    VOLUME = {155},
      YEAR = {2001},
    NUMBER = {1},
     PAGES = {1--16},
      ISSN = {0022-4049,1873-1376},
   MRCLASS = {18E15 (03C60 16B70 16D90)},
  MRNUMBER = {1804324},
MRREVIEWER = {Mike\ Prest},
       DOI = {10.1016/S0022-4049(99)00107-3},
       URL = {https://doi.org/10.1016/S0022-4049(99)00107-3},
}

@article {gab62,
    AUTHOR = {Gabriel, Pierre},
     TITLE = {Des cat\'egories ab\'eliennes},
   JOURNAL = {Bull. Soc. Math. France},
  FJOURNAL = {Bulletin de la Soci\'et\'e{} Math\'ematique de France},
    VOLUME = {90},
      YEAR = {1962},
     PAGES = {323--448},
      ISSN = {0037-9484},
   MRCLASS = {18.20},
  MRNUMBER = {232821},
       URL = {http://www.numdam.org/item?id=BSMF_1962__90__323_0},
}

@book {pop73,
    AUTHOR = {Popescu, N.},
     TITLE = {Abelian categories with applications to rings and modules},
    SERIES = {London Mathematical Society Monographs},
    VOLUME = {No. 3},
 PUBLISHER = {Academic Press, London-New York},
      YEAR = {1973},
     PAGES = {xii+467},
}

@article {loc_coh_grothendieck_spectrum,
    AUTHOR = {Krause, Henning},
     TITLE = {The spectrum of a locally coherent category},
   JOURNAL = {J. Pure Appl. Algebra},
  FJOURNAL = {Journal of Pure and Applied Algebra},
    VOLUME = {114},
      YEAR = {1997},
    NUMBER = {3},
     PAGES = {259--271},
      ISSN = {0022-4049,1873-1376},
       DOI = {10.1016/S0022-4049(95)00172-7},
       URL = {https://doi.org/10.1016/S0022-4049(95)00172-7},
}

@book {pp,
    AUTHOR = {Prest, Mike},
     TITLE = {Purity, spectra and localisation},
    SERIES = {Encyclopedia of Mathematics and its Applications},
    VOLUME = {121},
 PUBLISHER = {Cambridge University Press, Cambridge},
      YEAR = {2009},
     PAGES = {xxviii+769},
      ISBN = {978-0-521-87308-6},
       DOI = {10.1017/CBO9781139644242},
       URL = {https://doi.org/10.1017/CBO9781139644242},
}

@incollection{isolated,
      title={The importance of being isolated}, 
      author={Scott Balchin and Juan Omar Gómez and Greg Stevenson},
      year={2026},
      eprint={2605.19617},
      archivePrefix={arXiv},
      primaryClass={math.AC},
      url={https://arxiv.org/abs/2605.19617}, 
      note = "arXiv preprint: \url{https://arxiv.org/abs/2605.19617}"
}

@book {dst19,
    AUTHOR = {Dickmann, Max and Schwartz, Niels and Tressl, Marcus},
     TITLE = {Spectral spaces},
    SERIES = {New Mathematical Monographs},
    VOLUME = {35},
 PUBLISHER = {Cambridge University Press, Cambridge},
      YEAR = {2019},
     PAGES = {xvii+633},
      ISBN = {978-1-107-14672-3},
   MRCLASS = {54-02 (06B30 06D50 06F30 18B35 54B30 54H10)},
  MRNUMBER = {3929704},
MRREVIEWER = {Tomasz\ Kubiak},
       DOI = {10.1017/9781316543870},
       URL = {https://doi.org/10.1017/9781316543870},
}

@article {bhs23,
    AUTHOR = {Barthel, Tobias and Heard, Drew and Sanders, Beren},
     TITLE = {Stratification in tensor triangular geometry with applications
              to spectral {M}ackey functors},
   JOURNAL = {Camb. J. Math.},
  FJOURNAL = {Cambridge Journal of Mathematics},
    VOLUME = {11},
      YEAR = {2023},
    NUMBER = {4},
     PAGES = {829--915},
      ISSN = {2168-0930,2168-0949},
   MRCLASS = {18G80 (14F08 18F99 55P42 55P91 55U35)},
  MRNUMBER = {4650265},
       DOI = {10.4310/cjm.2023.v11.n4.a2},
       URL = {https://doi.org/10.4310/cjm.2023.v11.n4.a2},
}

@article {puninskii_superdecomposable,
    AUTHOR = {Puninski{\u\i}, G. E.},
     TITLE = {Superdecomposable purely injective modules over commutative
              valuation rings},
   JOURNAL = {Algebra i Logika},
  FJOURNAL = {Sibirski\u i\ Fond Algebry i Logiki. Rossi\u iskaya Akademiya
              Nauk. Sibirskoe Otdelenie. Institut Matematiki. Algebra i
              Logika},
    VOLUME = {31},
      YEAR = {1992},
    NUMBER = {6},
     PAGES = {655--671, 678},
      ISSN = {0373-9252},
       DOI = {10.1007/BF02261731},
       URL = {https://doi.org/10.1007/BF02261731},
}

@article {kru32,
    AUTHOR = {Krull, W.},
     TITLE = {Allgemeine {B}ewertungstheorie},
   JOURNAL = {J. Reine Angew. Math.},
  FJOURNAL = {Journal f\"{u}r die Reine und Angewandte Mathematik. [Crelle's
              Journal]},
    VOLUME = {167},
      YEAR = {1932},
     PAGES = {160--196},
      ISSN = {0075-4102,1435-5345},
       DOI = {10.1515/crll.1932.167.160},
       URL = {https://doi.org/10.1515/crll.1932.167.160},
}

@article {stevenson_absolutelyflatrings,
    AUTHOR = {Stevenson, G.},
     TITLE = {Derived categories of absolutely flat rings},
   JOURNAL = {Homology Homotopy Appl.},
  FJOURNAL = {Homology, Homotopy and Applications},
    VOLUME = {16},
      YEAR = {2014},
    NUMBER = {2},
     PAGES = {45--64},
      ISSN = {1532-0073},
   MRCLASS = {13D09},
  MRNUMBER = {3234500},
MRREVIEWER = {Sunil K. Chebolu},
       DOI = {10.4310/HHA.2014.v16.n2.a3},
       URL = {https://doi.org/10.4310/HHA.2014.v16.n2.a3},
}

@inproceedings {balmer-ICM,
    AUTHOR = {Balmer, Paul},
     TITLE = {Tensor triangular geometry},
 BOOKTITLE = {Proceedings of the {I}nternational {C}ongress of
              {M}athematicians. {V}olume {II}},
     PAGES = {85--112},
 PUBLISHER = {Hindustan Book Agency, New Delhi},
      YEAR = {2010},
      ISBN = {978-81-85931-08-3; 978-981-4324-32-8; 981-4324-32-9},
   MRCLASS = {18E30 (14F05 19G12 55P42)},
  MRNUMBER = {2827786},
MRREVIEWER = {Sunil\ K.\ Chebolu},
}

@article {neemanchromtower,
    AUTHOR = {Neeman, A.},
     TITLE = {The chromatic tower for {$D(R)$}},
      NOTE = {With an appendix by Marcel B{\"o}kstedt},
   JOURNAL = {Topology},
  FJOURNAL = {Topology. An International Journal of Mathematics},
    VOLUME = {31},
      YEAR = {1992},
    NUMBER = {3},
     PAGES = {519--532},
      ISSN = {0040-9383},
     CODEN = {TPLGAF},
   MRCLASS = {18E30},
  MRNUMBER = {1174255},
MRREVIEWER = {Steven E. Landsburg},
       DOI = {10.1016/0040-9383(92)90047-L},
       URL = {http://dx.doi.org/10.1016/0040-9383(92)90047-L},
}

@article {neemanBrownRep,
    AUTHOR = {Neeman, Amnon},
     TITLE = {The {G}rothendieck duality theorem via {B}ousfield's
              techniques and {B}rown representability},
   JOURNAL = {J. Amer. Math. Soc.},
  FJOURNAL = {Journal of the American Mathematical Society},
    VOLUME = {9},
      YEAR = {1996},
    NUMBER = {1},
     PAGES = {205--236},
      ISSN = {0894-0347,1088-6834},
   MRCLASS = {18E30 (14F05)},
  MRNUMBER = {1308405},
MRREVIEWER = {Luca\ Barbieri Viale},
       DOI = {10.1090/S0894-0347-96-00174-9},
       URL = {https://doi.org/10.1090/S0894-0347-96-00174-9},
}

@article {neemanOddball,
    AUTHOR = {Neeman, Amnon},
     TITLE = {Oddball {B}ousfield classes},
   JOURNAL = {Topology},
  FJOURNAL = {Topology. An International Journal of Mathematics},
    VOLUME = {39},
      YEAR = {2000},
    NUMBER = {5},
     PAGES = {931--935},
      ISSN = {0040-9383},
   MRCLASS = {18E30 (55P60)},
  MRNUMBER = {1763956},
MRREVIEWER = {Neil\ P.\ Strickland},
       DOI = {10.1016/S0040-9383(99)00040-3},
       URL = {https://doi.org/10.1016/S0040-9383(99)00040-3},
}

@article {bik11,
    AUTHOR = {Benson, D. and Iyengar, S. B. and Krause, H.},
     TITLE = {Stratifying triangulated categories},
   JOURNAL = {J. Topol.},
  FJOURNAL = {Journal of Topology},
    VOLUME = {4},
      YEAR = {2011},
    NUMBER = {3},
     PAGES = {641--666},
      ISSN = {1753-8416},
   MRCLASS = {18E30 (13D45 18G99 20J06 55P42)},
  MRNUMBER = {2832572 (2012m:18016)},
MRREVIEWER = {Guodong Zhou},
       DOI = {10.1112/jtopol/jtr017},
       URL = {http://dx.doi.org/10.1112/jtopol/jtr017},
}

@article {balmerfavi_idempotents,
    AUTHOR = {Balmer, P. and Favi, G.},
     TITLE = {Generalized tensor idempotents and the telescope conjecture},
   JOURNAL = {Proc. Lond. Math. Soc. (3)},
  FJOURNAL = {Proceedings of the London Mathematical Society. Third Series},
    VOLUME = {102},
      YEAR = {2011},
    NUMBER = {6},
     PAGES = {1161--1185},
      ISSN = {0024-6115},
   MRCLASS = {18E30 (14F05 55P60)},
  MRNUMBER = {2806103},
MRREVIEWER = {Fernando Muro},
       DOI = {10.1112/plms/pdq050},
       URL = {https://doi.org/10.1112/plms/pdq050},
}

@article {balmer_spectrum,
    AUTHOR = {Balmer, P.},
     TITLE = {The spectrum of prime ideals in tensor triangulated
              categories},
   JOURNAL = {J. Reine Angew. Math.},
  FJOURNAL = {Journal f\"ur die Reine und Angewandte Mathematik. [Crelle's
              Journal]},
    VOLUME = {588},
      YEAR = {2005},
     PAGES = {149--168},
      ISSN = {0075-4102},
     CODEN = {JRMAA8},
   MRCLASS = {18E30 (55P99)},
  MRNUMBER = {2196732},
MRREVIEWER = {Amnon Neeman},
       DOI = {10.1515/crll.2005.2005.588.149},
       URL = {http://dx.doi.org/10.1515/crll.2005.2005.588.149},
}

@article {hochster,
    AUTHOR = {Hochster, M.},
     TITLE = {Prime ideal structure in commutative rings},
   JOURNAL = {Trans. Amer. Math. Soc.},
  FJOURNAL = {Transactions of the American Mathematical Society},
    VOLUME = {142},
      YEAR = {1969},
     PAGES = {43--60},
      ISSN = {0002-9947},
   MRCLASS = {13.20 (14.00)},
  MRNUMBER = {251026},
MRREVIEWER = {S. S. Page},
       DOI = {10.2307/1995344},
       URL = {https://0-doi-org.pugwash.lib.warwick.ac.uk/10.2307/1995344},
}

@book {fs01,
    AUTHOR = {Fuchs, L. and Salce, L.},
     TITLE = {Modules over non-{N}oetherian domains},
    SERIES = {Mathematical Surveys and Monographs},
    VOLUME = {84},
 PUBLISHER = {American Mathematical Society, Providence, RI},
      YEAR = {2001},
     PAGES = {xvi+613},
      ISBN = {0-8218-1963-1},
       DOI = {10.1090/surv/084},
       URL = {https://doi.org/10.1090/surv/084},
}

@article {gldm1,
    AUTHOR = {Bazzoni, S. and {\v{S}}{t}\hspace{-0.125em}{'}\hspace{-0.125em}ov\'{\i}\v{c}ek, J.},
     TITLE = {Smashing localizations of rings of weak global dimension at
              most one},
   JOURNAL = {Adv. Math.},
  FJOURNAL = {Advances in Mathematics},
    VOLUME = {305},
      YEAR = {2017},
     PAGES = {351--401},
      ISSN = {0001-8708,1090-2082},
       DOI = {10.1016/j.aim.2016.09.028},
       URL = {https://doi.org/10.1016/j.aim.2016.09.028},
}

@article {keller94,
    AUTHOR = {Keller, B.},
     TITLE = {A remark on the generalized smashing conjecture},
   JOURNAL = {Manuscripta Math.},
  FJOURNAL = {Manuscripta Mathematica},
    VOLUME = {84},
      YEAR = {1994},
    NUMBER = {2},
     PAGES = {193--198},
      ISSN = {0025-2611,1432-1785},
       DOI = {10.1007/BF02567453},
       URL = {https://doi.org/10.1007/BF02567453},
}

@article {bt25,
    AUTHOR = {Balchin, Scott and Tecklenburg, Florian},
     TITLE = {Classifying smashing ideals in derived categories of valuation
              domains},
   JOURNAL = {J. Pure Appl. Algebra},
  FJOURNAL = {Journal of Pure and Applied Algebra},
    VOLUME = {229},
      YEAR = {2025},
    NUMBER = {3},
     PAGES = {Paper No. 107917, 26},
      ISSN = {0022-4049,1873-1376},
       DOI = {10.1016/j.jpaa.2025.107917},
       URL = {https://doi.org/10.1016/j.jpaa.2025.107917},
}

@article {bks20,
    AUTHOR = {Balmer, P. and Krause, H. and Stevenson, G.},
     TITLE = {The frame of smashing tensor-ideals},
   JOURNAL = {Math. Proc. Cambridge Philos. Soc.},
  FJOURNAL = {Mathematical Proceedings of the Cambridge Philosophical
              Society},
    VOLUME = {168},
      YEAR = {2020},
    NUMBER = {2},
     PAGES = {323--343},
      ISSN = {0305-0041,1469-8064},
       DOI = {10.1017/s0305004118000725},
       URL = {https://doi.org/10.1017/s0305004118000725},
}

@book {Glaz:1989,
    AUTHOR = {Glaz, Sarah},
     TITLE = {Commutative coherent rings},
    SERIES = {Lecture Notes in Mathematics},
    VOLUME = {1371},
 PUBLISHER = {Springer-Verlag, Berlin},
      YEAR = {1989},
     PAGES = {xii+347},
      ISBN = {3-540-51115-6},
   MRCLASS = {13-02 (13D05 13E99 13F05)},
  MRNUMBER = {999133},
MRREVIEWER = {Liam\ O'Carroll},
       DOI = {10.1007/BFb0084570},
       URL = {https://doi.org/10.1007/BFb0084570},
}

@article {Silver:1967,
    AUTHOR = {Silver, L.},
     TITLE = {Noncommutative localizations and applications},
   JOURNAL = {J. Algebra},
  FJOURNAL = {Journal of Algebra},
    VOLUME = {7},
      YEAR = {1967},
     PAGES = {44--76},
      ISSN = {0021-8693},
   MRCLASS = {16.40},
  MRNUMBER = {217114},
MRREVIEWER = {E.\ A.\ Walker},
       DOI = {10.1016/0021-8693(67)90067-1},
       URL = {https://doi.org/10.1016/0021-8693(67)90067-1},
}

@article {Lazard:1969,
    AUTHOR = {Lazard, Daniel},
     TITLE = {Autour de la platitude},
   JOURNAL = {Bull. Soc. Math. France},
  FJOURNAL = {Bulletin de la Soci\'et\'e{} Math\'ematique de France},
    VOLUME = {97},
      YEAR = {1969},
     PAGES = {81--128},
      ISSN = {0037-9484},
   MRCLASS = {16.90 (18.00)},
  MRNUMBER = {254100},
MRREVIEWER = {J.-P.\ Lafon},
       URL = {http://www.numdam.org/item?id=BSMF_1969__97__81_0},
}

@article {Thomason:1997,
    AUTHOR = {Thomason, R. W.},
     TITLE = {The classification of triangulated subcategories},
   JOURNAL = {Compositio Math.},
  FJOURNAL = {Compositio Mathematica},
    VOLUME = {105},
      YEAR = {1997},
    NUMBER = {1},
     PAGES = {1--27},
      ISSN = {0010-437X,1570-5846},
   MRCLASS = {18E30 (18F30)},
  MRNUMBER = {1436741},
MRREVIEWER = {Steven\ E.\ Landsburg},
       DOI = {10.1023/A:1017932514274},
       URL = {https://doi.org/10.1023/A:1017932514274},
}

@article {Buan-Krause-Solberg:support-var,
    AUTHOR = {Buan, Aslak Bakke and Krause, Henning and Solberg, \O yvind},
     TITLE = {Support varieties: an ideal approach},
   JOURNAL = {Homology Homotopy Appl.},
  FJOURNAL = {Homology, Homotopy and Applications},
    VOLUME = {9},
      YEAR = {2007},
    NUMBER = {1},
     PAGES = {45--74},
      ISSN = {1532-0073,1532-0081},
   MRCLASS = {18E30},
  MRNUMBER = {2280286},
       DOI = {10.4310/hha.2007.v9.n1.a2},
       URL = {https://doi.org/10.4310/hha.2007.v9.n1.a2},
}

@article {HPS:1997,
    AUTHOR = {Hovey, Mark and Palmieri, John H. and Strickland, Neil P.},
     TITLE = {Axiomatic stable homotopy theory},
   JOURNAL = {Mem. Amer. Math. Soc.},
  FJOURNAL = {Memoirs of the American Mathematical Society},
    VOLUME = {128},
      YEAR = {1997},
    NUMBER = {610},
     PAGES = {x+114},
      ISSN = {0065-9266,1947-6221},
   MRCLASS = {55U35 (18E30 18G35 55N20 55P42 55P60 55U15)},
  MRNUMBER = {1388895},
MRREVIEWER = {J.\ P. C. Greenlees},
       DOI = {10.1090/memo/0610},
       URL = {https://doi.org/10.1090/memo/0610},
}

@article {DG:2002,
    AUTHOR = {Dwyer, W. G. and Greenlees, J. P. C.},
     TITLE = {Complete modules and torsion modules},
   JOURNAL = {Amer. J. Math.},
  FJOURNAL = {American Journal of Mathematics},
    VOLUME = {124},
      YEAR = {2002},
    NUMBER = {1},
     PAGES = {199--220},
      ISSN = {0002-9327,1080-6377},
   MRCLASS = {16E30 (16D90 18E30)},
  MRNUMBER = {1879003},
MRREVIEWER = {Henning\ Krause},
       URL =
              {http://muse.jhu.edu/journals/american_journal_of_mathematics/v124/124.1dwyer.pdf},
}

@article {Herzog:1997,
    AUTHOR = {Herzog, Ivo},
     TITLE = {The {Z}iegler spectrum of a locally coherent {G}rothendieck
              category},
   JOURNAL = {Proc. London Math. Soc. (3)},
  FJOURNAL = {Proceedings of the London Mathematical Society. Third Series},
    VOLUME = {74},
      YEAR = {1997},
    NUMBER = {3},
     PAGES = {503--558},
      ISSN = {0024-6115,1460-244X},
   MRCLASS = {18E15 (03C60 16D90)},
  MRNUMBER = {1434441},
MRREVIEWER = {Henning\ Krause},
       DOI = {10.1112/S002461159700018X},
       URL = {https://doi.org/10.1112/S002461159700018X},
}
\end{document}